\documentclass[a4paper]{amsart}
\usepackage{mathrsfs}
\usepackage{bbm}
\usepackage{latexsym, amssymb,amsmath,amsthm}
\usepackage[all, knot]{xy}
\xyoption{arc}
\textheight=53pc \def\baselinestretch{1}
\renewcommand{\baselinestretch}{1.2} \topmargin=0mm
\newtheorem{thm}{Theorem}[section]
\newtheorem{cor}[thm]{Corollary}
\newtheorem{lem}[thm]{Lemma}
\newtheorem{prop}[thm]{Proposition}
\newtheorem{defn}[thm]{Definition}
\newtheorem{eg}[thm]{Example}

\newtheorem{rem}[thm]{Remark}

\def\<{\langle}
\def\>{\rangle}

\def\b{\beta}
\def\c{\cdot}

\def\D{\Delta}

\def\r{\rho}
\def\lr{\longrightarrow}

\def\o{\otimes}

\def\t{\tau}

\def\v{\varepsilon}

\def\HM{{}_{H}\mathscr{M}}
\def\BM{\mathrm{BM}}
\def\BP{\mathrm{BP}}
\def\Br{\mathrm{Br}}
\def\Gal{\mathrm{Gal}}
\def\Aut{\mathrm{Aut}}
\def\baseC{\mathscr{C}}

\def\autoequivtri{\Aut^{br}(\mathscr{Z}(\mathscr{C}),\mathscr{C})}
\date{}
\begin{document}
\renewcommand{\baselinestretch}{1.2}
\renewcommand{\arraystretch}{1.0}

\title{\bf Braided autoequivalences and quantum commutative bi-Galois objects }
\date{}
\author{Yinhuo Zhang}
\address{Department WNI, University of Hasselt, Universitaire Campus, 3590 Diepeenbeek, Belgium}
\email{yinhuo.zhang@uhasselt.be}
\author{Haixing Zhu}
\address{School of Economics and Management, Nanjing Forest University, Longpan Road 159, 210037, China}
\email{zhuhaixing@163.com}

\begin{abstract}
Let $(H,R)$ be a quasitriangular weak Hopf algebra over a field $k$. We show that there is a braided monoidal equivalence between the Yetter-Drinfeld module category $^H_H\mathscr{YD}$ over $H$ and  the  category of comodules over some braided Hopf algebra ${}_RH$ in the category $_H\mathscr{M}$. Based on this equivalence, we prove that every braided bi-Galois object $A$ over the braided Hopf algebra ${}_RH$ defines a braided autoequivalence of the  category $^H_H\mathscr{YD}$  if and only if $A$ is quantum commutative. In case $H$ is semisimple over an algebraically closed field, i.e. the fusion case, then every braided autoequivalence of $^H_H\mathscr{YD}$ trivializable on $_H\mathscr{M}$ is determined by such a quantum commutative Galois object.  The  quantum commutative Galois objects in $_H\mathscr{M}$ form a group measuring the Brauer group of $(H,R)$ as studied in \cite{Zh04} in the Hopf algebra case.
\end{abstract}

\maketitle

\section*{\bf{Introduction}}
Let $\mathscr{C}$ be a  braided fusion category  $\mathscr{C}$, that is,
a fusion category equipped with a braiding. Denote by $\mathscr{Z}(\mathscr{C})$ the \emph{Drinfeld center} of $\mathscr{C}$. The braided autoequivalences of $\mathscr{Z}(\mathscr{C})$ play important
roles in the study of braided fusion categories, see \cite{DN12,DGNO10,ENO10}.
For example, auto-equivalences  were used  to classify $G$-extensions of a given fusion category, see \cite{ENO10}. In order to
classify $G$-extensions of a given fusion category $\mathscr{C}$ using the
classical homotopy theory,  P. Etingof, D. Nikshych and V. Ostrik
introduced in \cite{ENO10} a 3-groupoid $\underline{\BP}(\mathscr{C})$, called the Brauer-Picard groupoid  of $\mathscr{C}$. This 3-groupoid can be truncated in the usual way into the Brauer-Picard group $\BP(\mathscr{C})$ of
$\mathscr{C}$, i.e. the group of the equivalence classes of invertible $\mathscr{C}$-bimodule categories. It turns out that there is a natural
 group isomorphism \cite[Thm 1.1]{ENO10}:
$$
\BP(\mathscr{C}) \cong  \Aut^{br}(\mathscr{Z}(\mathscr{C})),
$$
where  $\Aut^{br}(\mathscr{Z}(\mathscr{C}))$ is the group of isomorphism classes of braided autoequivalences of $\mathscr{Z}(\mathscr{C})$.
The name "Brauer-Picard group" speaks for itself that the group $\BP(\mathscr{C})$ has a close relation with the Brauer group $Br(\mathscr{C})$ of the category $\mathscr{C}$ which classifies the Azumaya algebras in $\mathscr{C}$, see \cite{VZ98}. In fact, every Azumaya algebra in $\mathscr{C}$ defines an invertible $\mathscr{C}$-bimodule category, so that $\Br(\mathscr{C})$ forms a subgroup of $\BP(\mathscr{C})$. The characterization of the Brauer group $\Br(\mathscr{C})$ in the group $\Aut^{br}(\mathscr{Z}(\mathscr{C}))$ has been done by A. Davydov and D. Nikshych in \cite{DN12}, where the braided autoequivalences corresponding to the Azumaya algebras are those trivializable on the base category $\mathscr{C}$, that is, $\Br(\baseC)\cong \Aut^{br}(\mathscr{Z}(\baseC),\baseC)$. 

Now we look at braided fusion categories  from the angle of weak
Hopf algebras. Let $k$ be an algebraically closed field. It is well known that a braided fusion category $\baseC$ is equivalent to the category $_H\mathscr{M}^{fd}$ of finite dimensional modules over some finite dimensional quasitriangular semisimple weak Hopf algebra $(H,R)$ over $k$, see
\cite{ENO05,NTV03,O03}. When the weak Hopf algebra $H$ happens to be
a Hopf algebra, we know that the Brauer group of $\baseC$ is the
Brauer group $\BM(H,R)$ of $(H,R)$ consisting of Azumaya $H$-module
algebras, see \cite{VZ98}.  In this case, the Brauer group  $\BM(H,R)$ can be characterized by the quantum commutative Galois objects over the braided Hopf algebra $_RH$,
 the transmutation of the quasitriangular Hopf algebra $(H,R)$, see \cite{Zh04}. In fact, we have the following general exact sequence
of groups:
$$
1 \longrightarrow \Br(k) \longrightarrow \BM(H,R)
\longrightarrow \Gal^{qc}(_RH), 
$$
where $\Gal^{qc} ({}_RH)$ is the group of quantum commutative
bi-Galois objects over ${}_RH$, and $k$ does not need to be algebraically closed. Now the question is whether the
group $\Gal^{qc}(_RH)$ is isomorphic to $\autoequivtri$, where $\baseC={}_H\mathscr{M}$. The answer
is positive, see \cite{DZ13}. The proof is based on the fact that an
autoequivalence of the comodule category over a Hopf algebra $H$ is
defined by a bi-Galois object over $H$. We don't know whether this
fact still holds for a weak Hopf algebra. However, one direction is
always true, that is, a bi-Galois object over a weak Hopf algebra
$H$ defines an autoequivalence of the comodule category over $H$.
In case $H$ is semisimple over an algebraically closed field, i.e. the braided category $_H\mathscr{M}^{fd}$ is a fusion category, we can show that both groups
$\Gal^{qc}(_RH)$ and $\autoequivtri$ are isomorphic to the Brauer
group $\BM(H,R)$, see \cite{ZZ2}. To obtain the isomorphisms, we
first construct a braided Hopf algebra $_RH$ from a quasitriangular weak Hopf algebra $(H,R)$. Unlike the Hopf algebra case, the original algebra $H$ can not be deformed into a Hopf algebra in the category of $H$-modules using Majid's transmutation theory. Here our braided Hopf algebra ${}_RH$ is nested on some centralizer subalgebra of $H$,  see \cite{LZ12}.

The next step is to use the braided Hopf algebra ${}_RH$ to describe the
Drinfeld center of the category of left $H$-modules using the category
of left ${}_RH$-comodules. Our result is the following (see Theorem 2.5).

\textbf{Theorem 1} \emph{Let $(H, R)$  be a quasitriangular weak Hopf algebra over a field $k$. Then the category  of Yetter-Drinfeld modules over $H$ is
equivalent to the category of  left  comodules over the braided Hopf algebra ${}_RH$ as a braided monoidal category.}

Following \cite[Thm 5.2]{Sc05} we know that a braided bi-Galois object $A$ over a braided Hopf algebra $\mathcal{H}$ in a braided monoidal category $\baseC$ defines an autoequivalence of the category $\baseC^\mathcal{H}$ of comodules over $\mathcal{H}$. Now we can apply this result to the braided Hopf algebra $_RH$ in the braided monoidal category $_H\mathscr{M}$ of a weak quasitriangular Hopf algebra $(H,R)$. Following Theorem 1, we know that the
category of  left comodules over ${}_RH$ is braided.  Thus a natural question arises: when is the autoequivalence defined by a braided bi-Galois object $A$ over $_RH$  a braided autoequivalence? Our answer is as follows ( see Theorem 3.6):

\textbf{Theorem 2} \emph{Let $(H, R )$  be a  quasi-triangular weak Hopf algebra over a field $k$. Assume that $A$ is a  braided bi-Galois object. Then the functor $A\Box-$ defines a braided autoequivalence of the category  of
Yetter-Drinfeld modules if and only if $A$ is quantum commutative.}

As a consequence, we obtain the following result:

\textbf{Theorem 3} \emph{Let $\mathscr{C}$ be a braided fusion category. Then the Drinfeld center of $\mathscr{C}$ is equivalent to the category of
finite dimensional left comodules over some braided Hopf algebra
${}_RH_\mathscr{C}$. If $A$ is a  braided bi-Galois object
over ${}_RH_\mathscr{C}$, then the functor $A\Box-$ defines a braided
autoequivalence of the Drinfeld center of  $\mathscr{C}$ trivializable on $\baseC$ if and only if $A$ is quantum commutative.}

The paper is organized as follows. In Section 1, we recall some
necessary definitions such as a weak Hopf algebra, a Yetter-Drinfeld
module and the Drinfeld center of a monoidal category.  In Section
2, we show that the category of Yetter-Drinfeld modules over a
quasitriangular weak Hopf algebra $(H,R)$ is equivalent to the
category of left comodules over the braided Hopf algebra $_RH$. In
Section 3, we show that a braided bi-Galois object $A$ over $_RH$
defines a braided autoequivalence of the category  of
Yetter-Drinfeld modules if and only if $A$ is quantum commutative. 
Such a braided autoequivalence is trivializable on the base category $_H\mathscr{M}$. In case $(H,R)$ is semisimple and $k$ is algebraically closed, then every braided auto-equivalence of $^H_H\mathscr{YD}$ trivializable on $_H\mathscr{M}$ is given by a quantum commutative Galois object over $_RH$. The proof will be given in the forthcoming paper \cite{ZZ2} as it is a consequence of the exact sequence of the Brauer group. In the last section, we compute the braided Hopf algebras $_RH$ of the face algebras defined by Hayashi in \cite{Ha99} and the quantum
commutative Galois objects over $_RH$.

\section{\bf{Preliminaries}}

Throughout this paper $k$ is a fixed field.
 Unless otherwise stated, unadorned tensor products will be over $k$. For a coalgebra over $k$, the coproduct will be denoted by $\D $. We adopt
 Sweedler's notation for the comultiplication  in \cite{Sw69}, e.g., $\D (a)=a_{1}\o a_{2}$.

 We assume that the reader is familiar with the notions of a
(braided) monoidal category, a ribbon or a  modular category (see
 \cite{Ka95})  as well as a braided fusion category  in \cite{ENO05}.
Moreover, we make free use of the notions of algebras, bialgebras
and Hopf algebras in a braided monoidal category, see \cite{Ma95}.

\subsection{Weak Hopf algebras}
We first recall the notion of a weak Hopf algebra. For more detail
on weak Hopf algebras, the reader is referred to  \cite{BNS99}. A
\emph{weak Hopf algebra} $H$ is a $k$-algebra $(H, m, \mu )$
 and a $k$-coalgebra $(H, \D , \varepsilon )$
 such that the following axioms hold:
\begin{enumerate}
\item[(i)] \quad $\Delta(hk)=\Delta(h)\Delta(k)$,
\item[(ii)] \quad $\Delta^{2}(1)=1_{1}\o 1_{2}1_{1^{'}}\o
 1_{2^{'}}=1_{1}\o 1_{1^{'}}1_{2}\o 1_{2^{'}}$,
\item[(iii)] \quad $\varepsilon(hkl)=\varepsilon(hk_{1})\varepsilon(k_{2}l)
 =\varepsilon(hk_{2})\varepsilon(k_{1}l)$,
\item[(iv)] \quad \mbox {There exists a $k$-linear map $S: H\lr H$,
 called the \emph{antipode}, satisfying}\nonumber \\
 \quad \quad \quad $h_{1}S(h_{2})=\v (1_{1}h)1_{2}, \quad \quad
 S(h_{1})h_{2}=1_{1}\v (h1_{2}),
 \quad \quad S(h)=S(h_{1})h_{2}S(h_{3})$,
\end{enumerate} for all $h, k, l\in H$.
We have two idempotent linear maps $\varepsilon_{t},\varepsilon_{s}$:
$H\longrightarrow H$ defined respectively by
\begin{eqnarray*}
\varepsilon_{t}(h)=\varepsilon(1_{1}h)1_{2}, \quad \quad
 \varepsilon_{s}(h)=1_{1}\varepsilon(h1_{2}),
 \end{eqnarray*}
called the \emph{target map}  and the \emph{source map} respectively. Their images $H_{t}$ and $H_{s}$ are called the target space and the source space respectively.  In fact,  $H_{t}$ and $H_{s}$ are Frobenius-separable subalgebras of $H$. Moreover,  the following equations hold:
\begin{eqnarray}
&& h_{1}\o h_{2}S(h_{3})=1_{1}h\o 1_{2},\\
&& S(h_{1})h_{2}\o h_{3}=1_{1}\o h1_{2},\\
&&  h_{1}\o S(h_{2})h_{3}=h1_{1}\o S(1_{2}),
\end{eqnarray}
\begin{eqnarray}
&& h_{1}S(h_{2})\o h_{3}=S(1_{1})\o 1_{2}h,\\
&& \v (g\v _t (h))=\v (gh)=\v (\v _s (g)h),\\
&& y1_{1}\o S(1_{2})=1_{1}\o S(1_{2})y,\\
&& zS(1_{1})\o 1_{2}=S(1_{1})\o 1_{2}z,
\end{eqnarray}
for $g,h\in H, y\in H_s$ and $ z\in H_{t}$.

\begin{rem}
\begin{enumerate}
\item[(i)] A weak Hopf algebra $H$ is an ordinary Hopf
algebra if and only if $\Delta (1)=1\o 1$ if and only if $\v $ is a
homomorphism if and only if $H_t=H_s=k1_H$.
\item[(ii)] The antipode $S$ is an anti-algebra isomorphism between $H_{t} $ and $H_{s}$.
\item[(iii)]  A weak Hopf algebra $H$ is called \emph{regular} if
 $S^2 (x)=x$ for all $x\in H_tH_s$.
\item[(iv)] Every weak Hopf algebra can be obtained by twisting
the comultiplication of a regular weak Hopf algebra and keeping the same algebra structure, see \cite{Ni02}.
\end{enumerate}
\end{rem}

In what follows, a weak Hopf algebra always means a regular weak
Hopf algebra. We recall the definition of a quasitriangular weak Hopf algebra from \cite{BNS99,NTV03}.

\begin{defn}
  \emph{Let $H$ be a weak Hopf algebra with a bijective
antipode $S$.  A \emph{quasi-triangular weak Hopf algebra} is a
 pair $(H,R)$, where
 $$R=R^1\otimes
 R^2\in\Delta^{cop}(1)(H\otimes_{k}H)\Delta(1),$$ satisfies the
 following conditions:
\begin{eqnarray}
&&(id\otimes \Delta)R=R_{13}R_{12},\ \  \\
&&(\Delta\otimes id)R=R_{13}R_{23},\ \ \\
&& \Delta^{cop}(h)R=R\Delta(h),\ \ \ \
\end{eqnarray}
where $h\in H$, $R_{12}=R\otimes 1$, $R_{23}=1\otimes R$, etc. Moreover,
there exists an element $\overline{R}\in
\Delta(1)(H\otimes_{k}H)\Delta^{cop}(1)$ such that $R\overline{R}=\Delta^{op}(1)$ and $\overline{R}R=\Delta(1)$. Such an element  $R$ is often called an \emph{$R$-matrix.}  In particular, $(H,R)$ is
called a \emph{triangular weak Hopf algebra} if
$\overline{R}=R^2\otimes R^1.$}
\end{defn}

 For any $y\in H_s$ and $z\in H_t$, the  following equations hold:
\begin{eqnarray}
&&(1\o z)R=R(z\o 1),\quad \quad\quad~(y\o1)R=R(1\o y),\\
&&(z\o1)R=(1\o S(z))R,\quad \quad(1\o y)R=(S(y)\o 1)R,\\
&&R(y\o1)=R(1\o S(y)),\quad \quad R(1\o z)=R(S(z)\o1),\\
&&(\v_s\o id)(R)=\D(1),\quad\quad\quad\quad(id\o\v_s)(R)=(S\o id)\D^{cop}(1),\\
&&(\v_t\o id)(R)=\D^{cop}(1),\quad\quad\quad(id\o\v_t)(R)=(S\o
id)\D(1).
\end{eqnarray}

\subsection{Modules over weak Hopf algebras}
Let $H$ be a weak Hopf algebra. Denote by ${}_H\mathscr{M}$
the category of left   $H$-modules. Then ${}_H\mathscr{M}$ forms a monoidal category $({}_H\mathscr{M},  \o_{t}, H_{t}, a, l, r)$ as follows:
\begin{enumerate}
 \item[(i)] for any two objects $M $ and $N$ in  $
  {}_H\mathscr {M}$,  $$M\o_{t}N=\{\sum m_i\o n_i \in M\o N |\sum  \Delta(1)(m_i\o n_i)=\sum m_i\o n_i\}.$$  Clearly, $M\o_{t}N = \Delta(1)(M\o N) \subseteq M\o N;$
\item[(ii)] for any two objects $M $ and $N$ in  $
  {}_H\mathscr{M}  $, the $H$-module structure on $M\o_t N$ is as follows: \ \ $h\c (m\o_t n)= h_1\c m\o_t h_2\c n$  for all $h\in H$ and  $m\in M $ and
$n\in N $;
\item[(iii)] $H_t$ is the unit object with $H$-action $h\c z=\v_t(hz)$, where
$h\in H, z\in H_t$, and the  $k$-linear maps $l_M $, $r_M $ and
their inverses are given by
\begin{eqnarray*}
&&l_M (1_1 \c z \o 1_2\c m)=z\c m,\ \  l^{-1}_M ( m)=1_1\c 1_H \o
1_2\c m\\
&&r_M (1_1 \c m \o 1_2\c z)=S(z)\c m,\ \  r^{-1}_M ( m)=1_1\c m \o
1_2,
\end{eqnarray*}
 for any $z\in H_t$ and $m\in M$, where   $M $ is an object in  $
 {}_H\mathscr{M} $.
\end{enumerate}

If $(H,R)$ is a  quasi-triangular weak Hopf algebra, then  the category  $
 {}_H\mathscr{M} $ can be equipped with a braiding
$C$ as  follows \cite[Prop. 5.2]{NTV03}:
$$
C_{M,N}(m\o_t n)=R^2\c n\o_t  R^1\c m,\  \mbox {for all }m\in M\
\mbox {and } n\in N,
$$
where    $M $  and $N$ are any two objects in ${}_H\mathscr{M} $.

\subsection{ Yetter-Drinfeld modules and the Drinfeld center}
\begin{defn} \emph{Let $H$ be a  weak Hopf algebra. A left $H$-module $M$ is called  a left \emph{Yetter-Drinfeld module}
 if $(M, \r^L)$ is a left $H$-comodule such that the following two conditions:
\begin{eqnarray*}
&\mathrm{(i)}& \r ^L(m)=m_{[-1]}\o m_{[0]}\in H\o
 _{t} V,\\
&\mathrm{(ii)}&  (h\c m) _{[-1]} \o (h\c m) _{[0]}=h_1m _{[-1]}S(h_3) \o
h_2 \c m _{[0]},
\end{eqnarray*}
 are satisfied for all $h\in H$ and $m\in M$. For a Yetter-Drinfeld module $M$, we have the identity:
 \begin{eqnarray}
m_{[-1]} \o m _{[0]}= m _{[-1]}S(1_2) \o 1_1 \c m _{[0]}, \  \mathrm{for}\ m\in M.
\end{eqnarray}}
\end{defn}

Denote by  ${}^{H}_{H}\mathscr{YD}$ the category of  left Yetter-Drinfeld
 modules. A Yetter-Drinfeld morphism is both left $H$-linear and left
 $H$-colinear. If the antipode $S$ is bijective, then ${}^{H}_{H}\mathscr{YD} $ is a braided monoidal category with the braiding given by
\begin{eqnarray*}
C_{V,W}(v\o w)=v_{[-1]}\c w\o v_{[0]},
\end{eqnarray*}
where $v\in V \in {}^{H}_{H}\mathscr{YD}$ and $w\in  W\in
{}^{H}_{H}\mathscr{YD}$. In particular, if  $(H, R)$ is a
quasi-triangular weak Hopf algebra, then every left $H$-module $M$
is automatically a  left Yetter-Drinfeld module with the following left
coaction:
$$
 \r^L (m)= R^2 \o R^1\c m, \ \forall m\ \in M.
$$
It is easy to see that the category ${}_{H}\mathscr{M}$ is a braided
monoidal subcategory of  ${}^{H}_{H}\mathscr{YD}$.

\begin{defn}  \emph{Let $H$ be a weak Hopf algebra with a bijective antipode $S$. An algebra $A$ in $ {}^{H}_{H}\mathscr { YD}$  is called
\emph{quantum commutative} if the following equation:
\begin{eqnarray*}
xy=(x_{[-1]}\c y  )x_{[0]}
\end{eqnarray*}
holds for all $x,y \in A$.} \end{defn}

\begin{defn}  \emph{Let $H$ be a weak Hopf algebra with a bijective antipode.  The left\emph{ Drinfeld center} $\mathscr{Z}_l
({}_{H}\mathscr{M})$ of  the monoidal category ${}_{H}\mathscr{M}$  is
the category, whose objects are pairs $(U, \nu_{U,-})$, where $U$ is
an object of $ {}_{H}\mathscr{M}$ and $\nu_{U,-}$ is a natural
family of isomorphisms,  called \emph{half-braidings}:
$$
\nu_{U,V}:  U\o V  \longrightarrow V\o U, \ \forall \  V\in
{}_{H}\mathscr{M}
$$
 satisfying the Hexagon Axiom.  Similarly, one can define the right
Drinfeld center of ${}_{H}\mathscr{M}$.} \end{defn}

\begin{lem}  \emph{\cite[Thm 2.6]{CWY05}} Let $H$ be a weak Hopf
algebra with bijective antipode.  Then $\mathscr{Z}_l
({}_{H}\mathscr{M})$ is equivalent to ${}^{H}_{H}\mathscr{YD}$ as a
braided monoidal category. \end{lem}

\section{\bf{The Drinfeld center of a quasi-triangular weak Hopf algebra}}

Let $H$ be a  quasi-triangular weak Hopf algebra. In this section,
we show  that there is a braided monoidal equivalence between
the Drinfeld center of the category of left $H$-modules and   the
category of left comodules over some braided Hopf algebra.

Denote by $C_H(H_s) $ the centralizer subalgebra of $H_s$ in $H$.
Clearly,\  $ C_H(H_s) = \{1_1 h S(1_2)|  \ \forall\ h\in H\}. $ The
algebra $C_H(H_s) $ is a left $H$-module algebra with the
adjoint action: $ \ \ h\c x =h_1xS(h_2) $ for all $ h\in H$ and $x\in
C_H(H_s) . $

Now we need Majid's transmutation theory in the case of a
quasi-triangular weak Hopf algebra. Recall Theorem 3.11 from \cite{LZ12}.

\begin{lem}\label{lem2.1}  Let $(H, R)$ be a
quasi-triangular weak Hopf algebra. Then $C_H(H_s)$ is a Hopf
algebra in the braided monoidal category ${}_{H}\mathscr{M}$ with
the following structures:

\emph{(i)}   the multiplication $\overline{\mu}$ and the unit $
\overline{\eta}$ are defined by:
\begin{eqnarray*}
&& \overline{\mu}: C_H(H_s)\o_t C_H(H_s) \longrightarrow C_H(H_s),\
\ a\o_t b
 \longmapsto  (1_1\c  a) (1_2 \c b),\\
 && \ \ \ \ \quad\quad \quad   \ \overline{\eta}= Id_{H_t}: H_t  \longrightarrow C_H(H_s), \ \ x \longmapsto
 x.
\end{eqnarray*}

\emph{(ii)}  The comultiplication $\overline{\D} $ and the counit $
\overline{\v}$ are given by:
\begin{eqnarray*}
&&\overline{\D}: C_H(H_s) \longrightarrow C_H(H_s) \o_t C_H(H_s),\ \
x \longmapsto
x_1S(R^2) \o R^1\c x_2,   \\
&& \quad \quad \quad  \overline{\v}= \v_t:  C_H(H_s) \longrightarrow
H_t,\ \ \ \ x \longmapsto \v_t (x).
\end{eqnarray*}

\emph{(iii)} The antipode is $\overline{S}$ defined by
\begin{eqnarray*}
\overline{S}: C_H(H_s)\longrightarrow C_H(H_s),\ \ x\longmapsto R^2
R'^2S(R^1xS(R'^1)).
\end{eqnarray*}
Moreover, ${}_RH$ is  cocommutative cocentral in the sense of
\emph{\cite{Sc10}}. \end{lem}

A Hopf algebra in a braided monoidal category is usually called a braided Hopf algebra in case the category does not need to be mentioned. In the sequel, we shall call the Hopf algebra  $C_H(H_s)$ in  ${}_{H}\mathscr{M}$  a \emph{ braided Hopf algebra} and denote it by ${}_RH$.

\begin{defn} \emph{\cite{Ma95}} \emph{Let $H$ be a  quasitriangular
weak Hopf algebra. Let $M$ be a left $H$-module. We call $(M, \r^l)$
a \emph{left $_RH$-comodule} in the category $\HM$ if
 $(M, \r^l)$ is a left $_RH$-comodule such that $\r^l $ is left
$H$-linear, i.e.,
 $$\r^l( h\c m) = h_1\c  m_{(-1)} \o h_2\c  m_{(0)},\ \forall h\in H, \ m\in M.$$}
\end{defn}
Similarly, one can  define a right  $_RH$-comodule and an $_RH$-bicomodule in the category $\HM$. For convenience, in the sequel,  {\it a left (right, bi-) $_RH$-comodule in the category $\HM$  will be called a left (right, bi-) $_RH$-comodule for short}.

Let $(M, \r^l)$ and $(N, \r^l)$ be two left ${}_R H$-comdules.
The tensor product $M\o_t N$ is a left ${}_R H$-comdule with the following
comodule structure:
$$
h\c (m\o n)=h_1\c m\o h_2\c n,\ \ \ \r^l (m\o n) = (\overline{\mu}\o
1\o 1)(1 \o C \o 1)(\r^l \o \r^l) (m\o n),
$$
where  $m\in M$, $n\in N$, $h\in H$ and $C$ is the braiding  in
${}_{H}\mathscr{M}$.

Denote  by  ${} ^{_RH} (_{H}\mathscr{M})$ the category of left
$_RH$-comodules. Note that a morphism in  ${} ^{_RH} (_{H}\mathscr{M})$ is both left $H$-linear and left  $_RH$-colinear. It is easy to see that the category ${}^{_RH} (_{H}\mathscr{M})$ is a monoidal category with the unit object given by $H_t$.

Now we discuss  the relation between the category  ${} ^{_RH}
(_{H}\mathscr{M})$ and the category of left Yetter-Drinfeld $H$-modules.

\begin{lem}\label{lem2.3}  Let $H$ be a  quasitriangular weak Hopf algebra. If  $(M, \r^l  )$ is a left $_RH$-comodule,
then $M$ is a left Yetter-Drinfeld $H$-module with the following
$H$-comodule structure:
\begin{eqnarray*}
&& \r^L ( m)=  m_{(-1)}R^2\o  R^1 \c m_{(0)}\in H\o M,
\end{eqnarray*}
where $\r^l (m) =m_{(-1)} \o m_{(0)}$ for all $m\in M.$ \end{lem}
\proof For any $m\in M,$ we first have
$$
 1_1m_{(-1)}R^2\o  1_2R^1 \c m_{(0)} = m_{(-1)}1_1R^2\o  1_2R^1 \c
 m_{(0)}=m_{(-1)}R^2\o R^1 \c m_{(0)}.
$$
 So  $ \r^L( M)  \in H\o_t M$. Namely, $\r^L$ is well-defined.

Next we verify   that $(M, \r^L)$ is a left $H$-comodule. For the
coassociativity, we have:
\begin{eqnarray*}
(1\o \r^L)\r^L &=&(1\o \r^L)(m_{(-1)}R^2\o R^1 \c m_{(0)})\\
&=&m_{(-1)}R^2\o  (R^1 \c m_{(0)})_{(-1)}q^2  \o  q^1\c (R^1 \c m_{(0)})_{(0)}\\
&=&m_{(-1)}R^2\o  (R^1_1 \c m_{(0)_{(-1)}})q^2  \o  q^1\c (R^1_2 \c m_{(0)})\\
&=&m_{(-1)_1}S(r^2)R^2\o  (R^1_1r^1 \c m_{(-1)_2 })  q^2  \o  q^1  R^1_2 \c m_{(0)}\\
&\stackrel{(9)}{=}&m_{(-1)_1}S(r^2)p^2R^2\o  (p^1r^1 \c m_{(-1)_2 })  q^2  \o  q^1  R^1 \c m_{(0)}\\
&\stackrel{(8)}{=}&m_{(-1)_1}\v_s(r^2) R^2\o  (r^1 \c m_{(-1)_2 })  q^2  \o  q^1  R^1 \c m_{(0)}\\
&\stackrel{(14)}{=}&m_{(-1)_1}1_1 R^2\o  (S(1_2) \c m_{(-1)_2 })
q^2  \o  q^1  R^1 \c m_{(0)}\\
&=&m_{(-1)_1}1_1 R^2\o  (  m_{(-1)_2 }S^2(1_2))  q^2  \o  q^1  R^1 \c m_{(0)}\\
&=&m_{(-1)_1}1_1 R^2\o  (  m_{(-1)_2 }1_2)  q^2  \o  q^1  R^1 \c m_{(0)}\\
&=&m_{(-1)_1} R^2\o  m_{(-1)_2 }  q^2  \o  q^1  R^1 \c m_{(0)}\\
&=&(m_{(-1)_1}R^2_1\o (m_{(-1)_2}R^2_2 \o R^1 \c m_{(0)})\\
&=&(\D\o 1)(m_{(-1)}R^2\o R^1 \c m_{(0)})=(\D\o 1)\r^L (m).
\end{eqnarray*}
The counit axiom holds as well because we have:
\begin{eqnarray*}
(\v\o 1)\r^L (m) &=&\v (m_{(-1)}R^2) (R^1 \c m_{(0)})
\stackrel{(5)}{=}\v (m_{(-1)}\v_t (R^2)) (R^1 \c m_{(0)})\\
&\stackrel{(15)}{=}&\v (m_{(-1)}1_2) (S (1_1) \c m_{(0)})=\v (m_{(-1)}S(1_1)) (1_2 \c m_{(0)})\\
&=&\v ( m_{(-1)}1_1) (1_2 \c m_{(0)})=\v ( 1_1m_{(-1)}) (1_2 \c
m_{(0)})=m,
\end{eqnarray*}
where the last equality follows from the counit of a left
$_RH$-comodule, namely,
$$
l\circ(\v_t\o 1)( m_{(-1)}\o  m_{(0)})=\v_t (m_{(-1)}) \c m_{(0)}=m.
$$

Finally, the compatible condition holds since
\begin{eqnarray*}
 h_1( m_{(-1)}R^2)\o h_2\c [R^1 \c m_{(0)}]
&=&h_11_1 m_{(-1)}S(1_2)R^2\o  h_2R^1 \c m_{(0)}\\
&\stackrel{(3)}{=} & h_1 m_{(-1)}S(h_2)h_3R^2\o  h_4R^1 \c m_{(0)}\\
&\stackrel{(10)}{=}& h_1 m_{(-1)}S(h_2)R^2h_4\o  R^1h_3 \c m_{(0)}\\
&=&(h_1\c  m)_{(-1)}R^2 h_2\o  R^1\c (h_1\c  m)_{(0)}.
\end{eqnarray*}
for all $m\in M$ and $h\in H.$   \qed

The following lemma says that the converse of Lemma 2.3 is also
true.

\begin{lem} \label{lem2.4}
Let $H$ be a  quasitriangular weak Hopf algebra
with an antipode $S$. If $(N, \r^L  )$ is a left Yetter-Drinfeld
module,  then $N$ is a left $_RH$-comodule with the following
structure:
\begin{eqnarray*}
&&  \r^l (n) = n _{[-1]}S(R^2) \o R^1\c n _{[0]},
\end{eqnarray*}
where $\r^L (n)= n _{[-1]} \o n _{[0]} $ for all $n\in N$. \end{lem}
\proof First of all, we need to  check that $\r^l$ is well-defined.
For any $n\in N$,
\begin{eqnarray*}
 1_1 [n _{[-1]}S(R^2)S(1_2) \o R^1\c n _{[0]}
&=&1_1n _{[-1]}S(1_2R^2) \o R^1\c n _{[0]}\\
& \stackrel{(11)}{=}&1_1n _{[-1]}S(R^2) \o R^11_2\c n _{[0]}\\
&=&1_1n _{[-1]}S(R^2) \o R^1\c (1_2\c n _{[0]})\\
&=&n _{[-1]}S(R^2) \o R^1\c  n _{[0]};\\
\ \ \ 1_1\c  [n _{[-1]}S(R^2)] \o 1_2R^1\c n _{[0]}&=& [n _{[-1]}S(R^2)]S(1_1) \o 1_2R^1\c n _{[0]}\\
&=& [n _{[-1]}S(1_1R^2)] \o 1_2R^1\c n _{[0]}\\
&=& [n _{[-1]}S(R^2)] \o R^1\c n _{[0]}.
\end{eqnarray*}
So $\r^l (N)\subset {}_RH\o_t N$. The $H$-linearity of  the map
$\r^l$ follows from the equations below:
\begin{eqnarray*}
h_1\c  [n _{[-1]}S(R^2)] \o h_2R^1\c n _{[0]}
&=& h_1n _{[-1]}S(R^2)S(h_2) \o h_3R^1\c n_{[0]}\\
&=&h_1n _{[-1]}S(h_2R^2) \o h_3R^1\c n _{[0]}\\
&\stackrel{(10)}{=}&h_1n _{[-1]}S(R^2h_3) \o R^1h_2\c n _{[0]}\\
&=&(h_1n _{[-1]}S(h_3))S(R^2) \o R^1\c (h_2\c n _{[0]})\\
&=&(h\c n) _{[-1]}S(R^2) \o R^1\c (h\c n) _{[0]}=\r^l(h\c n),
\end{eqnarray*}
for all $h\in H$. Now we show that $(N, \r^l)$ is a left $_RH$-comodule. For any $n\in
N,$
\begin{eqnarray*}
(1\o \r^l) \r^l (n)
&=& n _{[-1]}S(R^2) \o (R^1\c n _{[0]})_{[-1]} S(r^2)\o r^1\c (R^1\c n _{[0]})_{[0]}\\
&=& n _{[-1]}S(R^2) \o R^1_1 n _{[0]_{[-1]}} S(R^1_3) S(r^2)\o r^1\c (R^1_2\c n_{[0]_{[0]}})\\
&=& n _{[-1]_1}S(R^2) \o R^1_1 n _{[-1]_{2}} S(r^2R^1_3) \o r^1R^1_2\c n _{[0]}\\
& \stackrel{(10)}{=}& n_{[-1]_1}S(R^2) \o R^1_1 n _{[-1]_{2}} S(R^1_2r^2) \o R^1_3r^1\c n _{[0]}\\
& \stackrel{(9)}{=}& n _{[-1]_1}S(R^2q^2) \o R^1\c  [n _{[-1]_{2}} S(r^2)] \o q^1r^1\c n _{[0]}\\
& \stackrel{(8)}{=}& [n _{[-1]_1}S(r^2_2)]S(R^2) \o R^1\c  [n _{[-1]_{2}} S(r^2_1)] \o r^1\c n _{[0]}\\
&=& [n _{[-1]}S(r^2)]_1S(R^2) \o R^1\c  [n _{[-1]} S(r^2)]_2 \o
r^1\c n _{[0]} \end{eqnarray*}
\begin{eqnarray*}
&=& \overline{\D}[n _{[-1]}S(r^2)]\o r^1\c n _{[0]}=(\overline{\D}\o
1) \r^l (n).
\end{eqnarray*}
Hence the coassociativity holds. Finally, we verify that $\v_t$ satisfies the counit axiom:
\begin{eqnarray*}
&& \v_t (n _{[-1]}S(R^2))\c  (R^1\c n _{[0]})\\
&=&(\v_t (n _{[-1]}S(R^2))R^1)\c n _{[0]}
 \stackrel{(9)}{=}(1_2R^1)\c n _{[0]}\v (1_1n _{[-1]}\v_t [S(R^2)])\\
&=&(1_2R^1)\c n _{[0]}\v (1_1n _{[-1]}S[\v_s (R^2)])\stackrel{(14)}{=}(1_2S(1'_2))\c n _{[0]}\v (1_1n _{[-1]}S(1'_1))\\
&=&(1_21'_1)\c n _{[0]}\v (1_1n _{[-1]}1'_2)=(1_21'_1)\c n _{[0]}\v (1_1n _{[-1]}S(1'_2))\\
&=&1_2\c n _{[0]}\v (1_1n _{[-1]}S(1_3))=n.
\end{eqnarray*}

Therefore,   $(N, \r^l)$ is a left $_RH$-comodule. \qed

Combining Lemma \ref{lem2.3} and Lemma \ref{lem2.4}, we obtain the following theorem.

 \begin{thm} \label{thm2.5}
 Let $(H, R )$  be a  quasitriangular weak Hopf algebra.   Then
there is a monoidal equivalence  $\mathcal{F}$ from the category ${}
^{_RH} (_{H}\mathscr{M})$ of left $_RH$-comodules to the category
${}^{H}_{H}\mathscr{YD}$ of left Yetter-Drinfeld modules:
$$
\mathcal{F} : \  {} ^{_RH} (_{H}\mathscr{M}) \longrightarrow
{}^{H}_{H}{\mathscr{Y}}{\mathscr{D}},  \ (M,\r^l)  \longmapsto  (M,  \r^L),
$$
where  $\r^L$ is defined in Lemma \ref{lem2.3}. The quasi-inverse of  $\mathcal{F}$ is
$$
\mathcal{G} : \ {}^{H}_{H}\mathscr{{YD}} \longrightarrow  {} ^{_RH}
(_{H}\mathscr{M}),  \ (N, \r^L) \longmapsto  (N, \r^l ),
$$
where $\r^l$ is defined in Lemma \ref{lem2.4}. \end{thm}

\proof  We show first that $ \mathcal{G}  \mathcal{F}(M)=M$ for any
object $M$ in ${} ^{_RH} (_{H}\mathscr{M})$.  It is enough to verify
that $\r^l (m)= m_{(-1)}\o m_{(0)}$ for all $m\in M.$  Indeed,
\begin{eqnarray*}
\r^l(m)&=&  m _{[-1]}S(R^2) \o R^1\c m
_{[0]}\\
&=&  m_{(-1)}r^2  S(R^2) \o R^1\c [r^1 \c m_{(0)}]=  m_{(-1)}r^2  S(R^2) \o (R^1r^1) \c m_{(0)}\\
&\stackrel{(8)}{=}&  m_{(-1)}\v_t(R^2) \o R^1 \c m_{(0)}\stackrel{(15)}{=} m_{(-1)}1_2 \o S(1_1) \c m_{(0)}\\
&=&  S^{-1} (1_2)\c m_{(-1)} \o S(1_1) \c m_{(0)}=  S  (1_2)\c m_{(-1)} \o S(1_1) \c m_{(0)}\\
&=&  1_1\c m_{(-1)} \o 1_2 \c m_{(0)}=   m_{(-1)} \o  m_{(0)}.
\end{eqnarray*}

Next we show that  $ \mathcal{F}\mathcal{G}(N)  =N$  for any object
of  ${}^{H}_{H}\mathscr{YD}$. For all $n\in N,$
\begin{eqnarray*}
\r^L (n)&=& n_{(-1)}R^2\o  R^1 \c n_{(0)}= n_{ [-1]} S(r^2) R^2\o  R^1 \c (r^1\c n_{ [0]})\\
&=& n_{ [-1]} S(r^2) R^2\o  (R^1r^1) \c  n_{ [0]}\stackrel{(8)}{=} n_{ [-1]} \v_s( R^2) \o  R^1 \c   n_{ [0]} \\
&\stackrel{(14)}{=}& n_{ [-1]} 1_1 \o  S(1_2) \c   n_{ [0]} =n_{ [-1]} S(1_2) \o  1_1 \c   n_{ [0]} \\
&=& 1'_1n_{ [-1]} S(1_2) \o  1_1 \c ( 1'_2\c  n_{ [0]}) =1_1n_{ [-1]} S(1_3) \o  1_2 \c  n_{ [0]} \\
&=& n_{ [-1]} \o   n_{ [0]}.
\end{eqnarray*}

Finally, we verify that the triple $(\mathcal{G}, Id, Id)$ is
monoidal. It is clear that $\mathcal{G}(H_t)=H_t.$   For any two
left Yetter-Drinfeld modules $U$ and $V$,  the left ${_RH}$-comodule
structure  on $\mathcal{G}(U)\o \mathcal{G}(V)$ is as follows:
\begin{eqnarray*}
&\ \ \ & (\mu\o 1\o 1)(1 \o C \o 1)(\r^L \o \r^L) (u\o v) \\
 &=& (\mu\o 1\o 1)(1 \o C \o 1) (u_{(-1)} \o u_{(0)}\o n_{(-1)}\o
 v_{(0)})\\
  &=& (\mu\o 1\o 1) (u_{(-1)} \o R^2\c v_{(-1)}\o R^1\c u_{(0)} \o
  v_{(0)})\\
  &=& u_{(-1)}(R^2\c v_{(-1)})\o R^1\c u_{(0)} \o
 v_{(0)},
\end{eqnarray*}
where $u\in U$ and $v\in V.$ Now we have
\begin{eqnarray*}
&&  u_{(-1)}(R^2\c v_{(-1)})\o R^1\c u_{(0)} \o
 v_{(0)}\\
&=&(u _{[-1]}S(p^2))R^2{}_1 (v _{[-1]}S(q^2)) S(R^2{}_2) \o R^1 \c (p^1\c u_{[0]}) \o  q^1\c v_{[0]}\\
&\stackrel{(8)}{=}&(u _{[-1]}S(p^2))r^2 (v _{[-1]}S(q^2)) S(R^2) \o (R^1r^1 p^1)\c u_{[0]}) \o  q^1\c v_{[0]}\\
&\stackrel{(8)}{=}&u _{[-1]}\v_s(r^2) (v _{[-1]}S(q^2)) S(R^2) \o (R^1r^1 )\c u_{[0]}) \o  q^1\c v_{[0]}\\
&\stackrel{(14)}{=}&u _{[-1]}S(1_2) (v _{[-1]}S(q^2)) S(R^2) \o (R^11_1 )\c u_{[0]}) \o  q^1\c v_{[0]}\\
&=&u _{[-1]}S(1_2) (v _{[-1]}S(q^2)) S(R^2) \o R^1\c (1_1 \c u_{[0]}) \o  q^1\c v_{[0]}\\
&=&u _{[-1]} (v _{[-1]}S(q^2)) S(R^2) \o R^1\c  u_{[0]} \o  q^1\c v_{[0]}\\
&=&(u _{[-1]}v _{[-1]})S(R^2q^2) \o R^1\c  u_{[0]} \o  q^1\c v_{[0]}\\
&\stackrel{(9)}{=}& (u _{[-1]}v _{[-1]})S(R^2) \o R^1\c  (u_{[0]} \o   v_{[0]})\\
&=&(u\o_t v) _{[-1]}S(R^2) \o R^1\c (u\o_t v) _{[0]}=\r^l(u\o_t v).
\end{eqnarray*}
Hence, $ \mathcal{G}(U\o V) =\mathcal{G}(U)\o \mathcal{G}(V). $ The
verification of the other axioms for a monoidal functor are obvious.  \qed

Since the category of Yetter-Dinfeld modules is braided, the
equivalence $\mathcal{G}$ in Theorem 2.5 induces a braiding in the category of
left ${}_RH$-comodules such that the equivalence becomes braided.

\begin{cor}\label{cor2.6}  Let $(H, R )$  be a  quasitriangular weak
Hopf algebra.  Then the category of left
$_RH$-comodules is a braided monoidal category with a braiding
$\widetilde{C}$ given by
\begin{equation}\label{newbraiding}
 \widetilde{C }(u\o v) = u _{(-1)}R^2\c  v \o R^1\c u _{(0)}, \ \forall u\in U, \forall v\in V,
\end{equation}
where $U$ and $V $ are any two left ${}_RH$-comodules. The inverse
of $\widetilde{C }$ is given by
\begin{eqnarray*}
&& \widetilde{C}^{-1} ( v\o u)= R^1\c u_{(0)}\o S^{-1}(u_{(-1)}R^2)
\c v.
\end{eqnarray*}
Moreover,  the functor $\mathcal{G}$ in Theorem 2.5 gives a braided
monoidal equivalence. \end{cor}

\proof  Consider the following commutative diagram of isomorphisms:
$$
\begin{picture}(250,60)
\thinlines \put(38,45){\shortstack{$\mathcal{G}(U)\o _t
\mathcal{G}(V)
  $}}
\put(110,47){\vector(1,0){100}}
\put(220,45){\shortstack{$\mathcal{G}(U\o_t V))$}}
\put(65,42){\vector(0,-1){30}}
 \put(70,26){\shortstack{${C}_{
\mathcal{G}(U), \mathcal{G}(V)}$}} \put(200,26){\shortstack{$
\mathcal{G}(C_{U,V}) $}} \put(38,0){\shortstack{$\mathcal{G}(V) \o_t
\mathcal{G}(U)$}} \put(210,2){\vector(-1,0){100}}
\put(220,0){\shortstack{$\mathcal{G}(V \o_t U),$}}
\put(247,42){\vector(0,-1){30}}
\end{picture}\\
$$
where the horizontal isomorphisms  are given by $Id:  \mathcal{G}(X)\o \mathcal{G}(Y)\cong \mathcal{G}(Y\o X)$. Thus, the braiding
$\widetilde{C }$ is just the composition $Id^{-1} \circ C_{U,V} \circ
Id$. In fact, we have
\begin{eqnarray*}
  \widetilde{C}_{U, V} (u\o v)
   &=&Id  \circ  C_{U,V} \circ Id (u\o v)\\
   &=&Id  \circ  C_{U,V}  (u\o v)\\
   &=&Id   (u_{[-1]}\c v\o u_{[0]})\\
   &=&  u_{(-1)}R^2\c v\o R^1\c u_{(0)}.
\end{eqnarray*}
Similarly, one can obtain the inverse of $\widetilde{C}$. \qed

By Lemma 1.6 and Corollary 2.6 we obtain the following corollary.

\begin{cor} Let $(H, R )$  be a  quasitriangular weak Hopf algebra.
Then the Drinfeld center $\mathscr{Z}_l(_{H}\mathscr{M})$ of left
$H$-modules is equivalent to the category ${} ^{_RH}
(_{H}\mathscr{M})$  of left $_RH$-comodules as a braided monoidal
category.
\end{cor}

 As a special case,  we have the following corollary on a quasitriangular Hopf algebra:

\begin{cor} \label{Bos} Let $(H, R )$  be a quasitriangular Hopf algebra.  Then
the Drinfeld center of left $H$-modules is equivalent to the
category of left $_RH$-comodules as a braided monoidal category.
\end{cor}

\begin{rem} \emph{(i) when  $H$ is a finite dimensional quasitriangular
Hopf algebra, the functor $\mathcal{G}$ was first proved in
\cite{Zh04} to have a right adjoint.}

\emph{(ii)  Let $H$ be a finite dimensional quasitriangular Hopf
algebra. Following \cite[Prop 4.1]{Maj93} the quantum double $D(H)$
is isomorphic to a semidirect product $A\rtimes H$, where $A=H^*$ is
a braided Hopf algebra. By Corollary \ref{Bos} we may choose $A$ as the dual braided Hopf algebra $(_RH)^*$. Thus we have the following equivalences of braided monoidal categories:
$$
{} _{(_RH)^*\rtimes H}\mathscr{M} \cong  {}_{D(H)}\mathscr{M} \cong ^{H}
_{H}\mathscr{YD}\cong  \mathscr{Z}_l(_{H}\mathscr{M}).
$$
In case $H$ is infinite dimensional, we  have neither the usual quantum double $D(H)$ nor the dual braided Hopf algebra $(_RH)^*$. But Corollary \ref{Bos} always holds for any (finite or
infinite dimensional) quasitriangular  Hopf algebra over any field
(or even over a commutative ring). In particular, the Drinfeld center
is naturally equivalent to the category of comodules over
$BU_q(\mathfrak{g})$ studied in \cite{Maj93}.}
\end{rem}

\section{\bf{Quantum commutative Galois objects}}
In this section we study (braided) Galois objects over the Braided Hopf algebra $_RH$ of a  finite dimensional quasitriangular weak Hopf algebra $(H, R)$. We shall construct  braided autoequivalences of the Drinfeld center of ${}_H \mathscr{M}$  from braided bi-Galois objects.
For the details about braided Galois objects over a braided Hopf
algebra one is referred to \cite{Sc05, Sc10}.

Let $(H, R)$ be a finite dimensional quasitriangular weak Hopf algebra.  An object $X$ in ${}_H \mathscr{M}$ is flat if tensoring with $X$ preserves equalizers. A flat object $X$ is called \emph{faithfully flat} if tensoring with $X$ reflects isomorphisms. It is not hard to see that ${}_RH$ is
flat in the category ${}_H \mathscr{M}$ since ${}_RH$ is finite and
has a dual object.

\begin{defn} \emph{\cite{Sc05}} \emph{An algebra $A$ in   ${}_H
\mathscr{M}$   is called a \emph{left ${}_RH$-comodule algebra} if
$A$ is a left ${}_RH$-comodule such that the left comodule map $\rho^l$
satifies:
 $$
\r^l(ab)= a_{(-1)}(R^2\c b_{(-1)}) \o  (R^1\c a_{(0)})b_{(0)},
$$
for all  $a,b\in A$, where   $\r^l (a)=a_{(-1)}\o a_{(0)}$ . Namely,
$\r^l$ is an algebra map in ${}_H \mathscr{M}$. }
\end{defn}
 Similarly, an
algebra $A$ in   ${}_H \mathscr{M}$ is called a \emph{right
${}_RH$-comodule algebra} if $A$  with a right ${}_RH$-coaction
$\r^r$ is a right ${}_RH$-comodule such that
 $$
\r^r(ab)= a_{(0)} (R^2\c b_{(0)}) \o  (R^1\c a_{(1)})b_{(1)}, $$
where $a,b\in A$ and  $\r^r (a)=a_{(0)}\o a_{(1)}$. An
\emph{${}_RH$-bicomodule algebra} is both a left and a right $_RH$-comodule algebra such that the left and the right coations commute.

Now let $A$  be a right ${}_RH$-comodule algebra. The subalgebra
$$A_0= \{a\in A| \r^r (a)=a \o_t 1 =1_1a\otimes 1_2 \}
$$
is called the \emph{coinvariant subalgebra}. Similarly, one can define the  coinvariant subalgebra of a left ${}_RH$-comodule algebra. An $_RH$-coinvariant subalgebra $A_0$ is said to be \emph{trivial} if $A_0=H_t$.

\begin{defn} \emph{ \cite[Defn 2.1]{Sc10}  Let $A$  be a right
${}_RH$-comodule algebra.  $A$ is called a \emph{right braided
${}_RH$-Galois object} if $ A$ is  faithfully flat and   the
morphism
$$
\b : A \o _t A \longrightarrow A \o_t {}_RH, \ a\o_t b \longmapsto
ab_{(0)} \o_t b_{(1)}
$$
is an isomorphism.  Similarly, one  can define a \emph{left braided
${}_RH$-Galois object} and a \emph{braided bi-Galois object}.
}\end{defn}

The coinvariant subalgebra $A_0$ of a  right ${}_RH$-Galois object $A$ is trivial. So is the coinvariant subalgebra
 of a  left ${}_RH$-Galois object $A$.  Moreover, it is not hard to see that $({}_R H, \t_{{}_R H, -})$ is an object in
the Drinfeld center   $ \mathscr{Z}_l({}_H\mathscr{M})$, where
$\t_{{}_R H, -}$ is a half-braiding
\begin{eqnarray*}
  && \t_{{}_R H, M} :  {}_R H\o  M  \longrightarrow   M \o {}_R H, \  h\o m\longmapsto  r^2 R^1\c m \o r^1 h
  R^2.
  \end{eqnarray*}
Since ${}_RH$ is cocommutative cocentral, for any left
${}_RH$-comodule  $ (M, \r^l)$,  by   \cite{Sc10} there exists a natural right comodule structure induced by the half-braiding $\t_{{}_R H,
M}:$
$$
\r^r= \t_{{}_R H, M} \circ \r^l: M \longrightarrow {}_RH \o M
\longrightarrow M \o {}_RH,
$$
so that $(M,\r^l, \r^r )$ becomes an ${}_RH$-bicomodule. By \cite{Sc10} we
call $M$ \emph{cocommutative} if the right ${}_RH$-comodule is
induced by the left ${}_RH$-comodule as above.

\begin{defn} \emph{A cocommutative braided bi-Galois object $A$
is called a \emph{quantum commutative Galois object }if $A$ is
quantum commutative as an algebra in  ${}^{H}_{H}\mathscr{YD}$.}
\end{defn}

By Theorem 2.5 and Corollary 2.6, a left Yetter-Drinfeld
module is an $_RH$-bicomodule in $_H\mathscr{M}$. Thus we can consider the contensor product $M \Box_{_RH} N$, or $M\Box N$ for convenience, for two left Yetter-Drinfeld modules $M$ and $N$:
$$M\Box N=\{m\o_t n\in M \o_t N| \r^r (m) \o_t n = m \o_t \r^l (n)\},$$
or precisely,
\begin{eqnarray}
M\Box N= \{m\o n \in M\o_t N| r^2 \c m _{[0]}  \o r^1m _{[-1]} \o
n= m\o n _{[-1]}S(R^2) \o R^1\c n
 _{[0]} \}.
\end{eqnarray}
If $A$ is a braided ${}_RH$-bi-Galois object,   by \cite{Sc05} we have  an
isomorphism:
$$
\xi: (A\Box M )\o_t  (A\Box N ) \cong A \Box (M\o_t N),
$$
given by $\xi ((a\o m ) \o (b\o n) ) = a   (R^2\c b) \o R^1\c m \o n$,
for all $a,b \in A$, $m\in M$ and $b\in N$. Following \cite{Sc10}
the cotensor functor $A\Box-$ is a monoidal autoequivalence of  ${}
^{_RH} (_{H}\mathscr{M})$.

\begin{lem}\label{lem3.4}  Let $(H, R )$  be a  finite dimensional
quasitriangular weak Hopf algebra. If $A$ is a quantum commutative
Galois object, then the functor $A\Box-$ is a braided
autoequivalence of  ${} ^{_RH} (_{H}\mathscr {M}).$ \end{lem}

\proof  Let $A$ be a quantum commutative Galois object. By Theorem
2.5 and \cite{Ka95} it suffices to verify that the following diagram
is commutative:
$$
\begin{picture}(250,60)
\thinlines \put(15,45){\shortstack{$(A\Box M)\o _t (A\Box N)
  $}}
\put(110,47){\vector(1,0){100}} \put(220,45){\shortstack{$A\Box
(M\o_t N))$}} \put(65,42){\vector(0,-1){30}}
 \put(70,26){\shortstack{$\widetilde{C }_{
A\Box M, A\Box N}$}} \put(200,26){\shortstack{$ A\Box \widetilde{C
}_{M,N} $}} \put(320,26){\shortstack{$ (*) $}}
\put(15,0){\shortstack{$(A\Box N) \o_t (A\Box M)$}}
\put(110,2){\vector(1,0){100}} \put(220,0){\shortstack{$A\Box (N\o_t
M)$}} \put(247,42){\vector(0,-1){30}}
\end{picture}\\
$$
Indeed, on the one hand, for any $a\o m \in A\Box M$ and $b\o
n\in A\Box N$, we have:
\begin{eqnarray*}
&& [ \xi \circ (\widetilde{C }_{
A\Box M, A\Box N})][(a\o m ) \o (b\o n) ]\\
&= &  \xi [(a\o m )_{(-1)} r^2\c (b\o n)
\o  r^1 \c (a\o m )_{(0)}]\\
&= &  \xi [a_{(-1)} r^2\c  (b\o n)
\o  r^1 \c (a_{(0)}\o m )]\\
&=  &  \xi [a_{(-1)_1} r^2_1\c b \o a_{(-1)_2} r^2_2\c n
\o  r^1_1 \c a_{(0)}\o r^1_2 \c m )]\\
&=  &  [a_{(-1)_1} r^2_1\c b] [ R^2r^1_1 \c a_{(0)}]\o R^1a_{(-1)_2} r^2_2\c n \o r^1_2 \c m \\
&=  &  [a _{[-1]_1}S(q^2_2)] r^2_1\c b] [ R^2r^1_1 \c [q^1\c a _{[0]}]]\o R^1[a _{[-1]_2}S(q^2_1)] r^2_2\c n \o r^1_2 \c m \\
&=  &  [a _{[-1]_1}[S(q^2)r^2]_1\c b] [ R^2r^1_1 q^1\c a _{[0]}]]\o R^1[a _{[-1]_2}[S(q^2)r^2]_2] \c n \o r^1_2 \c m \\
&\stackrel{(9)}{=}&  [a _{[-1]_1}[S(q^2)r^2p^2]_1\c b] [ R^2r^1 q^1\c a _{[0]}]]\o R^1[a _{[-1]_2}[S(q^2)r^2p^2]_2] \c n \o p^1 \c m \\
&\stackrel{(8)}{=}&  [a _{[-1]_1}[\v_s(r^2)p^2]_1\c b] [ R^2r^1 \c a _{[0]}]]\o R^1[a _{[-1]_2}[\v_s(r^2)p^2]_2] \c n \o p^1 \c m \\
&\stackrel{(14)}{=}&  [a _{[-1]_1}[1_1p^2]_1\c b] [ R^2S(1_2) \c a _{[0]}]]\o R^1[a _{[-1]_2}[1_1p^2]_2] \c n \o p^1 \c m \\
&=&  [a _{[-1]_1}p^2_1\c b] [ R^2S(1_2) \c a _{[0]}]]\o R^1[a _{[-1]_2}1_1p^2_2] \c n \o p^1 \c m \\
&\stackrel{(16)}{=} &  [a _{[-1]_1}p^2_1\c b] [ R^2\c a _{[0]}]\o
R^1[a _{[-1]_2}p^2_2] \c n \o p^1 \c m ,
\end{eqnarray*}
where   Corollary 2.6  and Lemma 2.4 were used in the first and
fifth equality, respectively. On the other hand, we have:
\begin{eqnarray*}
&& (1 \o \widetilde{C }) \circ \xi [(a\o m ) \o (b\o n)]
\end{eqnarray*}
\begin{eqnarray*}
 &= &  a
(r^2\c b) \o \widetilde{C } (r^1\c m
\o n )\\
&= &  a   (r^2\c b) \o (r^1\c m)_{(-1)}W^2 \c n \o
W^1\c (r^1\c m)_{(0)}\\
&= &  a   (r^2\c b) \o (r^1_1\c m_{(-1)})W^2 \c n \o
W^1r^1_2\c  m_{(0)}\\
&= &  a   (r^2\c b) \o r^1_1 m_{(-1)}S(r^1_2)W^2 \c n \o
W^1r^1_3\c  m_{(0)}\\
&= &  a   (r^2\c b) \o r^1_1 m_{[-1]}S(R^2)S(r^1_2)W^2 \c n \o
W^1r^1_3R^1\c  m_{[0]}\\
&= &  a   (r^2\c b) \o r^1_1 m_{[-1]}S(r^1_2R^2)W^2 \c n \o
W^1r^1_3R^1\c  m_{[0]}\\
&= &  a   (r^2\c b) \o r^1_1 m_{[-1]}S(r^1_3)S(R^2)W^2 \c n \o
W^1R^1r^1_2\c  m_{[0]}\\
&= &  a   (r^2\c b) \o r^1_1 m_{[-1]}S(r^1_3)  \c n \o
r^1_2\c m_{[0]}\\
&= &  a   (r^2\c b) \o r^1_1 m_{[-1]}S(1_2)S(r^1_3)  \c n \o
r^1_21_1\c m_{[0]}\\
&\stackrel{(14)}{=}&  a   (r^2\c b) \o r^1_1
m_{[-1]}S(R^2)p^2S(r^1_3) \c n \o
r^1_2p^1R^1\c m_{[0]}\\
&\stackrel{(18)}{=}&   (R^2\c a_{[0]})   (r^2\c b) \o
r^1_1R^1a_{[-1]}p^2S(r^1_3) \c n \o
r^1_2p^1\c m\\
&= &  [(R^2\c a_{[0]})_{[-1]}  \c  (r^2\c b)](R^2\c a_{[0]})_{[0]}
\o r^1_1R^1a_{[-1]}p^2S(r^1_3) \c n \o
r^1_2p^1\c m\\
&= &  [(R^2_1 a_{[-1]_2} S(R^2_3)r^2\c b)] (R^2_2\c a_{[0]})\o
r^1_1R^1a_{[-1]_1}p^2S(r^1_3) \c n \o
r^1_2p^1\c m\\
&\stackrel{(8)}{=}&   [(R^2 a_{[-1]_2} S(Q^2_2)r^2\c b)] (Q^2_1\c
a_{[0]})\o r^1_1Q^1R^1a_{[-1]_1}p^2S(r^1_3) \c n \o
r^1_2p^1\c m\\
&\stackrel{(10)}{=}&   [(a_{[-1]_1}R^2  S(Q^2_2)r^2\c b)] (Q^2_1\c
a_{[0]})\o r^1_1Q^1a_{[-1]_2}R^1p^2S(r^1_3) \c n \o
r^1_2p^1\c m\\
&\stackrel{(9)}{=}&   [(a_{[-1]_1}R^2  S(U^2)V^2r^2\c b)] (Q^2\c
a_{[0]})\o V^1U^1Q^1a_{[-1]_2}R^1p^2S(r^1_2) \c n \o
r^1_1p^1\c m\\
&\stackrel{(14)}{=}&   [(a_{[-1]_1}R^2  1_1r^2\c b)] (Q^2\c
a_{[0]})\o S(1_2)Q^1a_{[-1]_2}R^1p^2S(r^1_2) \c n \o
r^1_1p^1\c m\\
&\stackrel{(11)}{=}&   [(a_{[-1]_1}R^2  r^2\c b)] (Q^2S(1_2)\c
a_{[0]})\o Q^1a_{[-1]_2}1_1R^1p^2S(r^1_2) \c n \o
r^1_1p^1\c m\\
&\stackrel{(16)}{=}&  [(a_{[-1]_1}R^2 r^2\c b)] (Q^2\c a_{[0]})\o
Q^1a_{[-1]_1}R^1p^2S(r^1_2) \c n \o
r^1_1p^1\c m\\
&\stackrel{(9)}{=}&  [(a_{[-1]_1}R^2 \c b)] (Q^2\c a_{[0]})\o
Q^1a_{[-1]_1}R^1_1p^2S(R^1_3) \c n \o
R^1_2p^1\c m\\
&\stackrel{(10)}{=}&  [(a_{[-1]_1}R^2 \c b)] (Q^2\c a_{[0]})\o
Q^1a_{[-1]_1}p^2R^1_2S(R^1_3) \c n \o
p^1R^1_1\c m\\
&= &  [(a_{[-1]_1}R^2 \c b)] (Q^2\c a_{[0]})\o Q^1a_{[-1]_1}p^21_2
\c n \o p^11_1R^1\c m\\
 &= &  [(a_{[-1]_1}R^2 \c b)] (Q^2\c a_{[0]})\o
Q^1a_{[-1]_1}p^2 \c n \o p^1R^1\c m,
\end{eqnarray*}
where  Corollary \ref{cor2.6}  and Lemma \ref{lem2.4} were used in the second and
fifth equality, respectively; the twelfth and the thirteenth
equations stemmed from the compatible condition and the quantum
commutativity respectively. Thus $$ \xi \circ (\widetilde{C }_{
A\Box M, A\Box N}) =(1 \o \widetilde{C }) \circ \xi.$$

Therefore, $A\Box -$ is a braided autoequivalence of ${} ^{_RH} (_{H}\mathscr {M})$.  \qed

\begin{lem}\label{lem3.5} Let $(H, R )$  be a  finite dimensional
quasitriangular weak Hopf algebra. Assume that  $A$ is a braided
bi-Galois object. If the functor $A\Box-$ defines a braided
autoequivalence of   ${} ^{_RH} (_{H}\mathscr{M})$, then $A$ is
quantum commutative. \end{lem}

\proof Assume that  the functor $A\Box-$ defines a braided
autoequivalence. We  have  the commutative diagram $(*)$. Let
$M$ and $N$ be two left ${}_RH$-comodules.  Following the proof of Lemma
\ref{lem3.4} we obtain the following equation:
\begin{eqnarray}\label{eq3.2}
 && a_{(0)}   (r^2\c b) \o r^1_1 a_{(1)}p^2S(r^1_3)\c n  \o r^1_2p^1\c
 m  \nonumber\\
&=&[a_{(-1)_1} r^2_1\c b] [ R^2r^1_1 \c a_{(0)}]\o R^1a_{(-1)_2}
r^2_2\c n \o r^1_2 \c m,
\end{eqnarray}
for all $a\o m \in A\Box M$ and $b\o n\in A\Box N$. Now let $M={}_RH$. Since $a_{(0)}\o a_{(1)}, b_{(0)}\o b_{(1)}
\in A\Box {}_RH $, we may substitute them for the elements $a\o m$ and $b\o n$ in the above equation  and obtain the following equation:
\begin{eqnarray*}
&& a_{(0)}   (r^2\c b_{(0)}) \o (r^1\c a_{(1)})_{[-1]} \c b_{(1)} \o
(r^1\c a_{(1)})_{[0]}\\
&=& [a_{(-1)_1} r^2_1\c b_{(0)}] [ R^2r^1_1 \c a_{(0)}]\o
R^1a_{(-1)_2} r^2_2\c b_{(1)} \o r^1_2 \c a_{(1)}.
\end{eqnarray*}
 Now we apply the map $1\o \v_t \o \v_t$ to the foregoing equality and obtain the following:
\begin{eqnarray*}
&& [a_{(-1)_1} r^2_1\c b_{(0)}] [ R^2r^1_1 \c a_{(0)}]\o  \v_t (
R^1a_{(-1)_2}
r^2_2\c b_{(1)}) \v_t (r^1_2 \c a_{(1)})\\
&=& a_{(0)}   (r^2\c b_{(0)}) \o  \v_t [(r^1\c a_{(1)})_{[-1]} \c
b_{(1)} ] \v_t [(r^1\c a_{(1)})_{[0]}].
\end{eqnarray*}
Since $\v_t$ is an algebra map in the category  ${}_H \mathscr{M}$
and $A$ is a right ${}_RH$-comodule algebra, we have
\begin{eqnarray*}
[a_{(-1)} r^2\c b] [ r^1 \c a_{(0)}]=ab,
\end{eqnarray*}
 which is equivalent to $$ ab = (a_{[-1]} \c
b)a_{[0]}.$$ Thus $A$ is quantum commutative.

Now we show that $A$ is cocommutative. Namely, we need to verify
that the right coaction $\r^r$ on $A$ is induced by its left
coaction $\r^l$ and the half-braiding. Note that the regular left
$H$-module $H$ has an induced Yetter-Drinfeld module structure,
where the comodule structure is given by
$$
\r^L (h)=R^2 \o R^1 h:=h_{[-1]} \o h_{[0]}.
$$
 By Lemma \ref{lem2.4} we have a left ${}_RH$-comodule structure on $H$, where  $\r^l (h) = 1\o_t h$ for any $h\in H.$  Namely, $(H, \r^l)$ is a trivial left
$_RH$-comodule. Now consider $A\Box {}_RH$ and $A\Box H$. Note that
$1_A \o_t 1_H \in A\Box H$ and $a_{(0)} \o a_{(1)}  \in A\Box
{}_RH$. Using Equation \eqref{eq3.2} we easily get:
\begin{eqnarray*}
 &&
 a_{(0)}   (r^2\c 1) \o r^1_1 a_{(1)}p^2S(r^1_3)  \o r^1_2p^1\c
a_{(2)}\\
&=& [a_{(-1)_1} r^2_1\c 1] [ R^2r^1_1 \c a_{(0)}]\o R^1a_{(-1)_2}
r^2_2\o r^1_2 \c a_{(1)}.
\end{eqnarray*}
Now on the one hand, we have:
\begin{eqnarray*}
 && a_{(0)}   (r^2\c 1_A) \o r^1_1 a_{(1)}p^2S(r^1_3)  \o
r^1_2p^1\c a_{(2)}\\
&=& a_{(0)}   (\v_t (r^2)\c 1_A) \o r^1_1 a_{(1)}p^2S(r^1_3)  \o
r^1_2p^1\c a_{(2)}\\
&\stackrel{(15)}{=}&  a_{(0)}   (1_2'\c 1_A) \o S(1'_1)1_1
a_{(1)}p^2S(1_3)  \o
1_2p^1\c a_{(2)}\\
&=& 1_1'\c a_{(0)}  \o  1'_21_1 a_{(1)}p^2S(1_3)  \o
1_2p^1\c a_{(2)}\\
&=& a_{(0)}  \o 1_1 a_{(1)}p^2S(1_3)  \o
1_2p^1\c a_{(2)}\\
&=& a_{(0)}  \o a_{(1)}p^2S(1_2)  \o
1_1p^1\c a_{(2)}\\
&\stackrel{(11)}{=}& a_{(0)}  \o a_{(1)}p^2   \o
 p^1\c a_{(2)}.
\end{eqnarray*}
On the other hand, we have:
\begin{eqnarray*}
&& [a_{(-1)_1} r^2_1\c 1_A] [ R^2r^1_1 \c a_{(0)}]\o R^1a_{(-1)_2}
r^2_2\o r^1_2 \c a_{(1)}
\end{eqnarray*}
\begin{eqnarray*}
&=&  [\v_t (a_{(-1)_1} r^2_1)\c 1_A] [ R^2r^1_1 \c a_{(0)}]\o
R^1a_{(-1)_2}
r^2_2\o r^1_2 \c a_{(1)}\\
&=&  [\v_t (a_{(-1)_1} \v_t(r^2_1))\c 1_A] [ R^2r^1_1 \c a_{(0)}]\o
R^1a_{(-1)_2}
r^2_2\o r^1_2 \c a_{(1)}\\
&\stackrel{(4)}{=}&  [\v_t (a_{(-1)_1} S(1_1))\c 1_A] [ R^2r^1_1
\c a_{(0)}]\o R^1a_{(-1)_2}
1_2r^2\o r^1_2 \c a_{(1)}\\
&=&  [\v_t (a_{(-1)_1} )\c 1_A] [ R^2r^1_1 \c a_{(0)}]\o
R^1a_{(-1)_2}
r^2\o r^1_2 \c a_{(1)}\\
&=&  [1_1 \c 1_A] [ R^2r^1_1 \c a_{(0)}]\o R^11_2a_{(-1)}
r^2\o r^1_2 \c a_{(1)}\\
&=& S(1_1) R^2r^1_1 \c a_{(0)}]\o R^11_2a_{(-1)}
r^2\o r^1_2 \c a_{(1)}\\
&\stackrel{(11)}{=}& R^2r^1_1 \c a_{(0)}\o R^1a_{(-1)} r^2\o r^1_2
\c a_{(1)}.
\end{eqnarray*}
Thus, the following equation holds:
$$  a_{(0)}  \o a_{(1)}p^2   \o
 p^1\c a_{(2)}=R^2r^1_1 \c a_{(0)}\o R^1a_{(-1)} r^2\o r^1_2 \c a_{(1)}\in A \o   H \o  H. $$
Applying the map $ (1\o 1\o \v)$ to right side of the above equation, we obtain:
\begin{eqnarray*}
&&  (1\o 1\o \v)(R^2r^1_1 \c a_{(0)}\o R^1a_{(-1)} r^2\o r^1_2 \c a_{(1)})\\
&\stackrel{(5)}{=}&R^2r^1_1 \c a_{(0)}\o  [R^1a_{(-1)} r^2]  \v( \v_s (r^1_2) a_{(1)}  S(r^1_3) )\\
&\stackrel{(2)}{=}&R^2r^1_1 \c a_{(0)}\o  [R^1a_{(-1)} r^2]  \v(   a_{(1)}  S(r^1_2) )\\
&=&R^2r^1_1 \c a_{(0)}\o  [R^1a_{(-1)} r^2]  \v(   a_{(1)} \v_t( S(r^1_2) ))\\
&=&R^2r^1_1 \c a_{(0)}\o  [R^1a_{(-1)} r^2]  \v(   a_{(1)} S( \v_s(r^1_2) ))\\
&\stackrel{(3)}{=}&R^2r^11_1 \c a_{(0)}\o  [R^1a_{(-1)} r^2]  \v(   a_{(1)} S^2( 1_2) ))\\
&=&R^2r^11_1 \c a_{(0)}\o  [R^1a_{(-1)} r^2]  \v(   a_{(1)} S( 1_2) ))\\
&=&R^2r^11_1 \c a_{(0)}\o  [R^1a_{(-1)} r^2]  \v(  S( 1_2) a_{(1)}  ))\\
&=&R^2r^11_1 \c a_{(0)}\o  [R^1a_{(-1)} r^2]  \v(  1_2  a_{(1)}  )\\
&=&R^2r^1 S(\v_t ( a_{(1)} )) \c a_{(0)}\o   [R^1a_{(-1)} r^2] \\
&=&R^2r^1 \c a_{(0)}\o   [R^1a_{(-1)} r^2 ],
\end{eqnarray*}
where the counit of a right ${}_RH$-comodule $A$ was used in the
last equality.   Now we have
\begin{eqnarray*}
R^2r^1 \c a_{(0)}\o   [R^1a_{(-1)} r^2  ]&=&  (1\o 1\o \v)(a_{(0)}
\o a_{(1)}p^2 \o
 p^1\c a_{(2)})\\
 &=& a_{(0)}  \o a_{(1)}p^2 \v(
 p^1\c  a_{(2)})\\
 &=& a_{(0)}  \o a_{(1)}p^2 \v(
 \v_s(p^1_1) a_{(2)}  S(p_2^1)]\\
  &=& a_{(0)}  \o a_{(1)}p^2 \v(
  a_{(2)}  S(p^1)]\\
&=& a_{(0)}  \o a_{(1)}p^2 \v(
  a_{(2)}  S(\v_s (p^1))]\\
  &\stackrel{(14)}{=}& a_{(0)}  \o a_{(1)}1_2 \v(
  a_{(2)}  S(1_1)]\\
   &=& a_{(0)}  \o a_{(1)}1_2 \v(
   1_1a_{(2)} )\\
  &=& a_{(0)}  \o a_{(1)}  \v_t(
  a_{(2)})\\
  &=& a_{(0)}  \o a_{(1)},
\end{eqnarray*}
where  the counit on  ${}_RH$  was used in the last equality.  This
means that a right ${}_RH$-comodule structure on $A$ is indeed induced by
its left ${}_RH$-coaction. Therefore, $A$ is a quantum commutative
 Galois object. \qed

Summarizing the foregoing arguments, we obtain the main result of this section:

\begin{thm}\label{thm3.6} Let $(H, R )$  be a  finite dimensional
quasitriangular weak Hopf algebra.  Assume that  $A$ is a  braided
bi-Galois object.  Then the functor $A\Box-$ defines a braided
autoequivalence of the category ${}^H_{H}\mathscr{YD}$ of
Yetter-Drinfeld modules if and only if $A$ is quantum commutative.
\end{thm}

\proof   Assume that  $A$ is a   braided bi-Galois object. By
Lemma \ref{lem3.4} and Lemma \ref{lem3.5}, the functor $A\Box-$ defines a braided autoequivalence of ${} ^{_RH} (_{H}\mathscr{M})$ if and only if $A$
is quantum commutative. Since   ${}
^{_RH} (_{H}\mathscr{M}) \cong {}^H_{H}\mathscr{YD}$  as braided
monoidal categories,   the functor $A\Box-$ induces a braided
autoequivalence of ${}^H_{H}\mathscr{YD}$  if and only if $A$ is
quantum commutative. \qed

Recall that the Drinfeld center $\mathscr{Z}_l({}_H\mathscr{M})$ is tensor equivalent to the Yetter-Drinfeld module category $ {}^H_{H}\mathscr{YD}$. Thus the functor $A\Box-$ defines a braided autoequivalence of the Drinfeld center  if and only if $A$ is quantum commutative. This holds as well for any quasitriangular Hopf algebra.

In order to deal with the case of a braided fusion category, we need
to restrict ourself to the category of finite dimensional
representations.  Denote by ${}_H\mathscr{M}^{f.d}$ and ${}^H_H\mathscr{YD}^{f.d}$ the category of finite dimensional left $H$-modules and  the category of finite dimensional left Yetter-Drinfeld modules respectively. Then $\mathscr{Z}_l({}_H\mathscr{M}^{f.d}) \cong {}^H_H\mathscr{YD}^{f.d}$. Thus, Theorem \ref{thm3.6} applies to ${}_H\mathscr{M}^{f.d}$.


\begin{cor} Let $\mathscr{C}$ be a braided fusion
category. Then the Drinfeld center of $\mathscr{C}$ is equivalent to
the category of finite dimensional left comoduels
over some braided Hopf algebra ${}_RH_\mathscr{C}$. Moreover, if $A$
is a braided bi-Galois object over ${}_RH_\mathscr{C}$, then the
cotensor functor $A\Box-$ defines a braided autoequivalence of the
Drinfeld center of $\mathscr{C}$ if and only if $A$ is quantum
commutative. \end{cor}

\proof Suppose that  $\mathscr{C}$ is a braided fusion category.
By \cite{O03} there exists a semisimple connected weak Hopf algebra
$H_\mathscr{C}$ such that $\mathscr{C}$ is (tensor) equivalent to the
category ${}_{H_\mathscr{C}}\mathscr{M}^{f.d}$ of finite dimensional
left $H_\mathscr{C}$-modules. Similar to the proof of Corollary 2.6,
one can endow the category ${}_{H_\mathscr{C}}\mathscr{M}^{f.d}$
with a braiding $\Phi$ such that the equivalence between the two categories   preserves the braidings.
Following \cite[Prop 5.2]{NTV03} one can define a quasitriangular structure
$R$ on $H_\mathscr{C}$ so that the braiding $\Phi$ of ${}_{H_\mathscr{C}}\mathscr{M}^{f.d}$ is induced by the quaisi-triangular structure $R$ of ${}_{H_\mathscr{C}}$.  \qed

To end this section,  we show that the quantum commutative Galois
objects over $_RH$ form a subgroup of the group of braided bi-Galois
objects (see \cite{Sc05}). In the Hopf algebra case, this subgroup
was defined in \cite{Zh04}. In what follows, we fix a finite
dimensional quasitriangular weak Hopf algebra $(H,R)$. A Galois
object means a braided bi-Galois object over the braided Hopf
algebra $_RH$ in the category $_H\mathscr{M}$. It is easy to see
that $_RH\Box-$ defines the identity functor of
$^{_RH}(_H\mathscr{M})$. So $_RH$ is a quantum commutative Galois
object.

\begin{lem}\label{lem3.8}  If $A$ and $B$ are two quantum
commutative  Galois objects, so is $A\Box B$. \end{lem}

\proof Assume that $A$ and $B$ are quantum commutative Galois
objects. Then  $A\Box -$ and  $B\Box -$ are braided
autoequivalences. So is the composition $(A\Box B )  \Box - $. Thus
by Proposition 3.3  $A\Box B$ is quantum commutative. \qed

Let $A$ a bi-Galois object $A$.  One can define a braided bi-Galois
object $A^{-1}=: (_RH\o A)^{co{_RH}}\subset {}_RH\o A^{op}$  such
that $A\Box A^{-1} \cong {}_RH$ and $A^{-1}\Box A \cong {}_RH$. For
more detail on $A^{-1}$, one may refer to \cite{Sc05}.

\begin{lem}\label{lem3.9}   If $A$ is a quantum commutative
Galois object, so is $A^{-1}$. \end{lem}

\proof  Suppose that  $A$ is a quantum commutative Galois object.
The functor  $A\Box -$ is a braided autoequivalence functor. It is easy to
see that $A^{-1}\Box-$ gives the inverse of the functor  $A\Box -$. By Lemma \ref{lem3.5},  the Galois object $A^{-1}$ is quantum commutative. \qed

Denote by $Gal^{qc} ({}_RH)$ the set of isomorphism classes
of the quantum commutative  Galois objects. Let  $[A]$ denote the
isomorphism class of a quantum commutative Galois object $A$. By
Lemma \ref{lem3.8} and Lemma \ref{lem3.9} we obtain the following.

\begin{thm} The set $Gal^{qc} ({}_RH)$
forms a group. The multiplication is induced by the cotensor product $\Box$ over ${}_RH$, the identity is given by $[{}_RH]$ and
the inverse of an element $[A]$ is represented by $A^{-1}$. \end{thm}

It is well-known  that the category ${}_{H}\mathscr{M}$ is braided subcategory of the Yetter-Drinfeld module category ${}^H_{H}\mathscr{YD}$. If $M$ is a left $H$-module.  Then $M$ possesses  a left $H$-comodule structure:
$$
\r^L (m)=R^2 \o R^1 \c m:=m_{[-1]} \o m_{[0]},
$$
so that  $(M, \r^L)$ is a left Yetter-Drinfeld module. It follows from Lemma \ref{lem2.4} that the induced left ${}_RH$-comodule structure on $M$ is trivial, namely,  $\r^l (m) = 1\o_t m $ for all $m\in M.$    If $A$ is a braided bi-Galois object, then $A\Box M\cong M$. Thus the functor $A\Box-$
restricts to the identity functor on the category of left
$H$-modules.

Now we consider the image of the group  $Gal^{qc} ({}_RH)$ in the group $\Aut^{br}({}^H_{H}\mathscr{YD})$ of braided autoequivalences of the Yetter-Drinfeld module category.

\begin{defn} \emph{ \cite[Defn 2.1]{DN12}}\emph{ A braided
autoequivalence $F$ of ${}^H_{H}\mathscr{YD}$ is called
\emph{trivializable} on ${} _{H}\mathscr{M}$ if the restriction
$F|_{{} _{H}\mathscr{M}}$ is isomorphic to the identity functor as a
braided tensor functor.  } \end{defn}

Denote by $\Aut^{br}({}^H_{H}\mathscr{YD}, {} _{H}\mathscr{M} )$  the
group of isomorphism classes of braided autoequivalences  of
${}^H_{H}\mathscr{YD}$  trivializable on $ {} _{H}\mathscr{M}$.

\begin{cor} The group $Gal^{qc}({}_RH)$ is a subgroup of
 the group  $\Aut^{br}({}^H_{H}\mathscr{YD}, {}_{H}\mathscr{M})$.
 \end{cor}

We expect that the two groups are isomorphic for any finite
dimensional quasitriangular weak Hopf algebras $(H,R)$. This is the
case when $H$ is a Hopf algebra, see \cite{DZ13}. In case $H$ is
semisimple over an algebraically closed field, i.e. the fusion case,
the two groups are indeed isomorphic (to the Brauer group of the
braided fusion category), see \cite{ZZ2} or \cite{Zh12}.

\begin{eg}  \emph{Let $k$ be a field with $ch(k)\neq 2$. Let $H_4
$ be the Sweedler 4-dimensional Hopf algebra over $k$. Namely, $H_4$
is generated by two elements $g$ and $h$ satisfying
$$g^2=1, \ \ h^2=0,\ \ gh+hg=0.$$
The comultiplication, the counit and the antipode are given as
follows:
\begin{eqnarray*}
&& \D (g)= g\o g , \ \  \D (h)= 1\o h+h\o g\\
&& \v(g)=1,\ \ S(g)=g, \ \v(h)=0,\ \ S(h)=gh.
\end{eqnarray*}
It is known that $H_4$ has a quasitriangular structure $R_0$. All
quantum commutative Galois objects were computed in \cite {Zh04}.
Moreover, the group $Gal^{qc} (_{R_0}H)$ is isomorphic to $\Gamma
\rtimes Z_2$, where $\Gamma\cong k^+\times K^{\bullet }/K^{\bullet
2}.$ }\end{eg}

\section{\bf{Face algebras}}

In this section we compute the groups of quantum commutative Galois objects of a class of weak Hopf algebras, namely, the face algebras introduced by Hayashi in \cite{Ha99}.

Let $N\geq 2$ be an integer and $\mathbb{Z}_N$ the cyclic group
$\mathbb{Z}/N\mathbb{Z}.$ Let $\omega\in \mathbb{C}$ be a primitive $N^{th}$ root of unity. Let $H$  be the $\mathbb{C}$-linear span  of
$\{X^i_j (s) | i,j, s \in \mathbb{Z}_N \}$.  $H$ is  a
quasitriangular weak Hopf algebra equipped with the following
 structures:
\begin{eqnarray*}
&& \D (X^i_j (s))= \sum _{p+q=s} X^i_j (p)\o
X^{i+p}_{j+p} (q), \   \v(X^i_j (s))=\delta_{s,0},\\
&&X^i_j (p)X^k_l (q)=\delta_{j,k}\delta_{p,q}X^i_l (p),\ \   1=\sum_{i,p}X^i_i (p),\\
&& S(X^i_j (p))=X^{j+p}_{i+p} (-p),\\
&& R_1 \o R_2= \sum_{i,j,p}X^i_j (p) \o X^j_{j+p} (i-j)\omega^{-p(i-j)},\\
&&R'_1 \o R'_2= \sum_{i,j,p}  X^{j+p}_{i+p}(-p) \o X^{j}_{j+p}
(i-j)\omega^{-p(i-j)},
\end{eqnarray*}
 where  the target subalgebra $H_t$ of $H$ is  the $\mathbb{C}$-linear span  of $\{\sum_p
X^i_i(p)| i\in \mathbb{Z}_N\}$. Denote  by $1^i$ the sum $\sum _p
X^i_i (p)$ for all $i\in \mathbb{Z}_N.$ Then $H_t$
 is commutative and is equal to the direct sum $ \bigoplus_{i\in \mathbb{Z}_N } \mathbb{C}1^i$.

 Now we compute  the braided Hopf algebra ${}_RH $.

 \begin{lem}\label{lem4.1} The braided Hopf algebra ${}_RH $  is equal to the $\mathbb{C}$-linear span of $\{X^i_i (p)|i, p \in \mathbb{Z}_N\}$ equipped with the following structures:
\begin{eqnarray*}
&& \D' (X^k_k (s))= \sum_{w+q=s} X^k_k (w)\o
X^{k}_{k} (q),\ \  \v_t(X^i_i (s))=\delta_{s,0} \sum_p X^i_i (p),\\
&&X^i_i (p)X^k_k (q)=\delta_{i,k}\delta_{p,q} X^i_i (p) ,\
1=\sum_{i,p}X^i_i (p),\\
&& S (X^{k}_{k} (s)) =X^{k}_{k} (-s).
\end{eqnarray*}
 \end{lem}
\proof  Note that  $ \D(1_H) = \D(\sum_{i,s}X^i_i (s)) =
\sum_{i,s}\sum_{p+q=s} X^i_i (p) \o X^{i+p}_{i+p} (q).  $   We have
\begin{eqnarray*}
1_1 X^{m}_{n} (r) S(1_2)&=&  \sum_{i,s}\sum _{p+q=s}  X^i_i (p)
X^{m}_{n} (r)
 S( X^{i+p}_{i+p} (q) )\\
 &=&  \sum_{i,s}\sum_{p+q=s}  X^i_i (p)
X^{m}_{n} (r)
 X^{i+p+q}_{i+p+q} (-q)\\
 &=&  \sum_{i,s}\sum_{p+q=s}  \delta_{i,m}\delta_{n,i+p+q}  \delta_{p,r}\delta_{-q,r} X^i
_{i+p+q} (p)\\
&=&  \sum_{i} \delta_{i,m}\delta_{n,i} X^i _{i} (r)=\delta_{m,n} X^m
_{n} (r),
\end{eqnarray*}
for all $m,n,r\in \mathbb{Z}_N $.  So $_RH$ is the $\mathbb{C}$-linear
span of $\{X^i_i (p)| i, p \in V\}$.

Using the expression $ \D(R^1)\o R^2 = \sum_{i,j,p}\sum_{u+v=p} X^i_j(u) \o
X^{i+u}_{j+u}(v) \o X^j_{j+p}(i-j) \omega^{-p(i-j)}, $
  we compute the deformed comultiplication as follows:
\begin{eqnarray*}
&& \D' (X^k_k(s)) \\
&=& \sum_{w+q=s}   X^k_k(w)S(R^2) \o R^1 \c X^{k+w}_{k+w}(q)\\
&=& \sum_{w+q=s}   X^k_k(w)S(R^2) \o R^1_1
X^{k+w}_{k+w}(q)S(R^1_2)\\
&=& \sum_{w+q=s}\sum_{i,j,p}\sum_{u+v=p} X^k_k(w)S( X^j_{j+p}(i-j))
\o  X^i_j(u) X^{k+w}_{k+w}(q)S(X^{i+u}_{j+u}(v)) \omega^{-p(i-j)}\\
&=& \sum_{w+q=s}\sum_{i,j,p}\sum_{u+v=p} X^k_k(w) X_i^{i+p}(j-i)) \o
X^i_j(u)
X^{k+w}_{k+w}(q)X^{j+u+v}_{i+u+v}(-v)) \omega^{-p(i-j)} \\
&=& \sum_{w+q=s}\sum_{i,j,p}\sum_{u+v=p}\delta_{w,j-i}
\delta_{k,i+p} X^k_i(w) \o \delta_{u,q}\delta_{q,-v}
\delta_{j,k+w}\delta_{k+w,j+u+v}X^i_{i+u+v}(u)
 \omega^{-p(i-j)}\\
&=& \sum_{w+q=s}\sum_{i,j}\delta_{w,j-i} \delta_{k,i}
\delta_{j,k+w}X^k_i(w) \o X^i_{i}(q)\\
 &=& \sum_{w+q=s}\sum_{j}\delta_{w,j-k}
 \delta_{j,k+w}X^k_k(w) \o X^k_{k}\\
 &=&\sum_{w+q=s}
  X^k_k(w) \o X^k_{k}(q).
\end{eqnarray*}

By Lemma \ref{lem2.1} the antipode is given by\  $ \overline{S}(x) =
R^2R'^2S^2(R'^1)S(R^1x). $ For convenience,  we first compute
$R'^2S^2(R'^1)$. Indeed,
\begin{eqnarray*} R'^2S^2(R'^1)&=&
\sum_{i,j,p} X^j_{j+p}(i-j) S^2 (X^i_j (p)) \omega^{-p(i-j)}\\
&=&  \sum_{i,j,p} X^j_{j+p}(i-j) X^i_j (p) \omega^{-p(i-j)}\\
&=& \sum_{i,j,p} \delta_{j+p,i}  X^j_{j} (i-j)\omega^{-p(i-j)}.
\end{eqnarray*}
Now we have
\begin{eqnarray*}
\overline{S}(X^{k}_{k} (-s))&=&R^2 R'^2S(R^1X^{k}_{k} (s)S(R'^1))\\
&=& \sum_{i,j,p } \sum_{i',j',p' } \delta_{j'+p',i'}X^j_{j+p} (i-j)
X^{j'}_{j'} (i'-j')
 S(X^i_j (p) X^k_k (s) )\omega^{-[p(i-j)+p'(i'-j')]}\\
 &=& \sum_{i,j,p } \sum_{i',j',p' } \delta_{j'+p',i'}\delta_{j,k}\delta_{p,s}X^j_{j+p} (i-j)
X^{j'}_{j'} (i'-j')
 S(X^i_k (s) )\omega^{-[p(i-j)+p'(i'-j')]}\\
 &=& \sum_{i } \sum_{i',j',p' } \delta_{j'+p',i'}X^k_{k+s} (i-k)
X^{j'}_{j'} (i'-j')
 S(X^i_k (s) )\omega^{-[s(i-k)+p'(i'-j')]}\\
 &=& \sum_{i } \sum_{i',j',p' } \delta_{j'+p',i'}X^k_{k+s} (i-k)
X^{j'}_{j'} (i'-j')
 X^{k+s}_{i+s} (-s) \omega^{-[s(i-k)+p'(i'-j')]}\\
  &=& \sum_{i } \sum_{i',j',p' } \delta_{j'+p',i'}  \delta_{i-k,i'-j'}  \delta_{i'-j',-s}\delta_{k+s,j'}
 X^{k}_{i+s} (-s) \omega^{-[s(i-k)+p'(i'-j')]}\\
 \end{eqnarray*}
\begin{eqnarray*}
 &=& \sum_{i } \sum_{j',p' }    \delta_{i-k,p'}  \delta_{p',-s}\delta_{k+s,j'}
 X^{k}_{i+s} (-s) \omega^{-[s(i-k)+p'p']}\\
 &=& \sum_{i } \sum_{j' }    \delta_{i-k,-s} \delta_{k+s,j'}
 X^{k}_{i+s} (-s) \omega^{-[s(i-k)+(-s)(-s)]}\\
 &=& \sum_{i }  \delta_{i-k,-s}
 X^{k}_{i+s} (-s) \omega^{-[s(i-k)+(-s)(-s)]}\\
 &=&
 X^{k}_{k} (-s) \omega^{-[s(-s)+(-s)(-s)]}=X^{k}_{k} (-s).
 \end{eqnarray*}
Thus, the proof is completed. \qed

Take $i\in \mathbb{Z}_N$. Define $H^i$ to be the $\mathbb{C}$-linear
span of $\{ X^i_i (p)| p\in \mathbb{Z}_N\}$. It is obvious that $H^i$
is a subalgebra of ${}_RH$ with unity $ 1^i$. Moreover, ${}_RH $ is
the direct sum of all these $H^i$, i.e., $ {}_RH = \bigoplus_{i\in
\mathbb{Z}_N} H^i. $ We will show that  every $H^i$ is also an
ordinary Hopf algebra and so ${}_RH $ is actually the direct sum of
all these Hopf algebras. In order to verify that every $H^i$ can be
equipped with a coalgebra structure, we need to decompose the vector
space ${}_RH\o_t {}_RH. $

\begin{lem}\label{lem4.2}
${}_RH\o_t {}_RH = \bigoplus_{i\in \mathbb{Z}_N} (H^i \o H^i). $
 \end{lem}

\proof It is equivalent to show that
$$1_1\c X^a_a (b) \o 1_2\c X^u_u(w)=  \delta_{u,a}   X^a_a
  (b)    \o
X^{u}_{u}  (w),$$ for all $a,b,u,w \in \mathbb{Z}_N$. Indeed, we have
 \begin{eqnarray*}
 1_1\c X^a_a (b) \o 1_2\c X^u_u(w)
&=& \sum_{i,s} \sum_{p+q=s}  X^i_i (p) \c X^a_a (b)    \o
X^{i+p}_{i+p} (q) \c X^u_u(w)\\
&=& \sum_{i,s} \sum_{p+q=s}  \delta_{i,a} \delta_{p,0}  X^i_i
  (b)    \o   \delta_{i+p,u} \delta_{q,0}
X^{i+p}_{i+p}  (w)\\
&=& \sum_{i}   \delta_{i,a}   X^i_i
  (b)    \o   \delta_{i,u}
X^{i}_{i}  (w)=  \delta_{u,a}   X^a_a
  (b)    \o
X^{u}_{u}  (w),
 \end{eqnarray*}
for all $a,b,u,w \in \mathbb{Z}_N$. \qed

\begin{lem}\label{lem4.3} For all $i\in \mathbb{Z}_N$,
 $H^i$ is a
coalgebra over $\mathbb{C} 1^i $ with the following structures:
\begin{eqnarray*}
&& \D' (X^i_i (s))= \sum_{w+q=s} X^i_i (w)\o X^{i}_{i} (q),
\\
  &&\v_t(X^i_i (s))=\delta_{s,0} \sum_p
X^i_i (p).
\end{eqnarray*}
\end{lem}

\proof Follows from Lemma \ref{lem4.1} and Lemma \ref{lem4.2}. \qed

\begin{prop}\label{prop4.4} For all $i\in \mathbb{Z}_N $,  $H^i$ is a commutative and cocommutative Hopf
algebra over $\mathbb{C} 1^i $ equipped  with  the following
structures:
\begin{eqnarray*}
&&X^i_i (p)X^i_i (q)= \delta_{p,q} X^i_i (p) ,\ \ \ \ \
1_{H^i}=1^i,\\
&& \D' (X^i_i (s))= \sum_{w+q=s} X^i_i (w)\o
X^{i}_{i} (q), \  \\
&&\v_t(X^i_i (s))=\delta_{s,0} \sum_p X^i_i (p), \ \  S (X^{i}_{i}
(s)) =X^{i}_{i} (-s).
\end{eqnarray*}

\end{prop}

\proof Since we know already  that $H^i$ is both an algebra and a
coalgebra, it remains to be proved that    $\D'$ and $\v_t$ are
multiplicative, and that the axioms of   the antipode $S$ hold. We
first check that $\D'$ is multiplicative. Indeed,
 \begin{eqnarray*}
\D'(X^i_i (s)) \D''(X^i_i (t)) &=& [\sum_{p+q=s} X^i_i (p)\o
X^{i}_{i} (q)] [\sum_{p'+q'=t} X^i_i (p')\o
X^{i}_{i} (q')]\\
&=& \sum_{p+q=s}\sum_{p'+q'=t}
[X^i_i (p)X^i_i (p') \o X^{i}_{i} (q) X^{i}_{i} (q')]\\
&=& \sum_{p+q=s}\sum_{p'+q'=t} \delta_{p,p'}\delta_{q,q'}  [X^i_i
(p)
\o X^{i}_{i} (q)]\\
&=& \delta_{s,t} \sum_{p+q=s} X^i_i (p)\o  X^{i}_{i} (q) \\
&=&\D'(X^i_i (s) X^i_i (t)),
\end{eqnarray*}
for all $i,s,u,t\in \mathbb{Z}_N.$

Note that $\D' (1)=1 \o_t 1.$ It follows from Lemma \ref{lem4.2}  that
  $\D' (1^i)=   1^i \o 1^i$.

Next we verify that $\v_t$ is an algebra map.  For all $ s,t \in
\mathbb{Z}_N$, we have
\begin{eqnarray*}
\v_t (X^i_i (s)  ) \v_t (  X^i_i(t)) &= &
\delta_{s,0}\delta_{t,0}(\sum_p
X^i_i(p)) (\sum_q X^i_i(q))\\
&= & \delta_{s,0}\delta_{t,0} (\sum_p X^i_i(p)) =
\delta_{s,t}\delta_{s,0}\v_t (X^i_i (s))\\
&= &  \v_t (X^i_i (s) X^i_i(t)).
\end{eqnarray*}

Finally, we prove that the  antipode axioms  hold. Indeed,
\begin{eqnarray*}
m(1\o S) \D'' (X^i_i (s)) &=& \sum_{p+q=s} X^i_i(p)
S(X^i_i(q))= \sum_{p+q=s} X^i_i(p) X^i_i(-q)\\
&=& \delta_{p,-q}\sum_{p+q=s} X^i_i(p)= \delta_{s,0}\sum_{p\in
\mathbb{Z}_N} X^i_i(p)=\v_t (X^i_i (s)).
\end{eqnarray*}
for any $s\in  \mathbb{Z}_N.$ Similarly, we also have
\begin{eqnarray*}
 \sum_{w+q=s} S(X^i_i (w)X^{i}_{i} (q)&=&\sum_{w+q=s} X^i_i (-w)X^{i}_{i} (q))
 =\sum_{w+q=s} \delta_{-w,q} X^i_i (q)
  \end{eqnarray*}

\begin{eqnarray*}
&=&\sum_{q} \delta_{s,0} X^i_i (q)=\v_t (X^i_i (s)).
\end{eqnarray*}
Hence,  $H^i $ is an ordinary  Hopf algebra over $\mathbb{C} 1^i $.
\qed

In fact, $H^i$ is isomorphic to the dual Hopf algebra of the group Hopf algebra $k\mathbb{Z}_N$.

\begin{cor}\label{cor4.5}  The braided Hopf algebra $ {}_RH$ has a decomposition:
$$
_RH =\bigoplus_{i\in \mathbb{Z}_N} H^i,
$$
where $H^i$ is a Hopf algebra over $\mathbb{C} 1^i $ with unity
$1^i$.
 Moreover, there exists a Hopf algebra isomorphism from
$H^i$ to  $H^j$ defined by
$$
 \iota^j_i : X^i_i(p)\longmapsto X^j_j(p),
$$
for all $i,j,p\in \mathbb{Z}_N.$
\end{cor}

\proof Follows from  Proposition \ref{prop4.4}. \qed

Corollary \ref{cor4.5} indicates that braided bi-Galois objects over
${}_RH$ can be obtained from bi-Galois objects over a Hopf algebra
$H^i$.

Let the notations be as above.
Let $A$ be a quantum commutative Galois object over $_RH$.  Corollary \ref{cor4.5} implies that there is  a decomposition: $A =\bigoplus _{i\in
\mathbb{Z}_N} A^i,$ where $\r^r (A^i) \in A^i \o H^i$. Furthermore,
every $A^i$ is just a Galois object over $H^i$ (automatically a bi-Galois object as $H^i$ is cocommutative). Conversely, given a
Galois object $A'$ over  Hopf algebra $H^i$ for some $i\in
\mathbb{Z}_N$,  we can get a quantum commutative Galois object over
${}_R H$ as the direct sum $\bigoplus_{i\in \mathbb{Z}_N} A'^i$, where
every algebra $A'^i$ is a copy of $A'$.  Now we state the relation
between quantum commutative Galois object over ${}_R H$ and  Galois
object over $H^i$ as follows:

\begin{prop}
Let $A$ be a $\mathbb{C}$-algebra with unity. Then $A$ is   a
quantum commutative Galois object over $_RH$ if and only if $A $ is the direct
sum $\bigoplus_{i\in \mathbb{Z}_N}  A^i$, where every $A^i$ is an
$H^i$-Galois object. Moreover, there exists a group isomorphism
$$
\Omega:  Gal^{qc}({}_RH) \longrightarrow  Gal(H^i),\ \  A\longmapsto
A^i,
$$
for any fixed $i\in \mathbb{Z}_N$. The inverse of $\Omega$ is given as
follows: $$ \Omega':  Gal(H^i) \longrightarrow  Gal^{qc}({}_RH),\ \
A'\longmapsto \bigoplus_{i\in \mathbb{Z}_N} A'^i.
$$
\end{prop}

The detailed proof of the statement above is given in
\cite{Zh12} following a tedious and long computation. So the group
$Gal^{qc}({}_RH)$ can be obtained by computing the group $Gal(H^i)$ of Galois objects over  $H^i$. Since the Hopf algebra $H^i$ is commutative and cocommutative isomorphic to $k\mathbb{Z}_N$, we know that the group $Gal(H^i)$ is actually given by the second Galois cohomology group $H^2(\mathbb{Z}_N, k)$.

\section*{\bf{Acknowledgement}}

This work forms a part of the PhD thesis of the second named author at University of Hasselt. He would like to thank BOF of UHasselt for the financial support: BOF09-DOC009-R-1964.

 \end{document}